\documentclass[a4paper,review]{cas-sc}

\usepackage[section]{placeins} 
\usepackage{amssymb}
\usepackage{amsthm}
\usepackage{amsmath}
\usepackage{upgreek}
\usepackage{pdflscape}
\usepackage{listings}
\usepackage{multirow}
\usepackage[percent]{overpic}
\usepackage{multicol}
\usepackage{color}
\usepackage{stmaryrd}
\usepackage{xfrac}
\usepackage[capitalise]{cleveref}
\usepackage{booktabs}
\usepackage{tabularx}
\usepackage[section]{placeins} 
\usepackage[section]{algorithm}
\usepackage{calc}
\usepackage[titletoc]{appendix}
\usepackage{siunitx}
\usepackage{scalefnt}
\usepackage[sort&compress,square,numbers]{natbib}
\usepackage{cancel}

\usepackage{graphicx}
\usepackage{graphics}
\usepackage{wrapfig}
\usepackage{float}
\usepackage{subfig}
\usepackage[percent]{overpic}
\usepackage{varwidth}
\usepackage{tikz}
\usepackage{pgfplots}
\usepackage{adjustbox}
\usetikzlibrary{arrows,matrix,positioning,fit}
\usetikzlibrary{shapes,positioning}
\usetikzlibrary{backgrounds}
\usepackage{tikz-layers}
\usepgfplotslibrary{fillbetween}
\usetikzlibrary{intersections}
\pgfplotsset{compat=1.14}
\usepgfplotslibrary{colorbrewer}
\usepgfplotslibrary{patchplots}
\usepgfplotslibrary[colorbrewer]
\usetikzlibrary{pgfplots.colorbrewer}
\usetikzlibrary[pgfplots.colorbrewer]
\usepgfplotslibrary{units}
\usetikzlibrary{spy}
\usepackage{pgfplotstable}
\usepackage{arrayjobx}
\graphicspath{ {./figs/} }
\usetikzlibrary{external}

\newlength\myheight
\newlength\mydepth
\settototalheight\myheight{Xygp}
\newcommand*\inlinegraphics[1]{%
  \settototalheight\myheight{Xygp}%
  \settodepth\mydepth{Xygp}%
  \raisebox{-\mydepth}{\includegraphics[height=\myheight]{#1}}%
}
\newcommand\orcid[1]{\href{https://orcid.org/#1}{\inlinegraphics{orcid_16x16.png}}}

\makeatletter
\def\BState{\State\hskip-\ALG@thistlm}
\makeatother

\newdefinition{definition}{Definition}[section]

\def\lim{{\xi}}    

\colorlet{PlotColor1}{Spectral-A}
\colorlet{PlotColor3}{Spectral-L}
\colorlet{PlotColor2}{Spectral-O}
\colorlet{PlotColor4}{Spectral-D}
\colorlet{PlotColor5}{black}

\begin{document}


\title[mode=title]{A method for bounding high-order finite element functions: Applications to mesh validity and bounds-preserving limiters}
\shorttitle{Bounding high-order finite element functions}

\author[1]{Tarik Dzanic}[orcid=0000-0003-3791-1134]
\cormark[1]
\cortext[cor1]{Corresponding author}
\ead{dzanic1@llnl.gov}
\author[1]{Tzanio Kolev}[orcid=0000-0002-2810-3090]
\author[1]{Ketan Mittal}[orcid=0000-0002-2062-852X]
\shortauthors{T. Dzanic \textit{et al.}}

\address[1]{Center for Applied Scientific Computing, Lawrence Livermore National Lab, Livermore, CA 94550, USA}

\begin{abstract}
We introduce a novel method for bounding high-order multi-dimensional polynomials in finite element approximations. The method involves precomputing optimal piecewise-linear bounding boxes for polynomial basis functions, which can then be used to locally bound any combination of these basis functions. This approach can be applied to any element/basis type at any approximation order, can provide local (i.e., subcell) extremum bounds to a desired level of accuracy, and can be evaluated efficiently on-the-fly in simulations. Furthermore, we show that this approach generally yields more accurate bounds in comparison to traditional methods based on convex hull properties (e.g., Bernstein polynomials). The efficacy of this technique is shown in applications such as mesh validity checks and optimization for high-order curved meshes, where positivity of the element Jacobian determinant can be ensured throughout the entire element, and continuously bounds-preserving limiters for hyperbolic systems, which can enforce maximum principle bounds across the entire solution polynomial. 
\end{abstract}

\begin{keywords}
Finite element methods \sep
High-order \sep
Mesh validity \sep
Bounds-preserving \sep
Bounding box
\end{keywords}



\maketitle

\section{Introduction}
\label{sec:intro}
Bounding the extrema of polynomials remains an open problem in a variety of computational fields, including numerical analysis, scientific computing, and engineering design, where polynomials frequently appear in approximations of complex physical models and geometries. For example, in computer graphics, polynomial representations are used to visualize curves and surfaces (e.g., through B\'ezier curves and B-splines), where bounding their extrema is crucial for collision detection, rendering optimization, and surface smoothness analysis. Similarly, in numerical analysis and scientific computing, high-order polynomial approximations are used in finite element methods, where ensuring bounded behavior is critical for stability and accuracy in simulations. Despite their importance, determining tight bounds on polynomial extrema remains challenging, particularly in high-dimensional spaces or when dealing with polynomials of high degree. Traditional approaches such as interval arithmetic, convex hull techniques, and sum-of-squares optimization provide partial solutions but often suffer from computational complexity or overly conservative bounds.

This work investigates techniques for bounding polynomial extrema of high-order finite element-type approximations, although with broader applicability to problems in computer graphics and contact mechanics. In particular, we focus on two primary applications: i) mesh validity checks for high-order curved meshes; and ii) continuously bounds-preserving limiters for high-order finite element approximations of hyperbolic systems. In regard to the former, mesh validity checks are essential to ensure that high-order curved elements maintain geometric integrity, preventing issues such as element inversion. Unlike low-order meshes, where element validity is straightforward to assess, high-order elements introduce additional challenges due to their curved nature, which is represented by high-order polynomial approximations in terms of an element transformation matrix. Effective bounding techniques for the determinant of this transformation matrix help certify that these element transformations remain well-posed, which is necessary to maintaining accuracy and stability in finite element simulations for both static and deforming meshes.

For the latter, bounds-preserving limiters are crucial for maintaining the stability and physical consistency of high-order finite element approximations of hyperbolic systems. However, the vast majority of limiters, which enforce bounds discretely at nodal or quadrature points, do not guarantee that the polynomial representation remains bounded at arbitrary locations within an element. This becomes problematic in applications requiring solution evaluation at new points, such as adaptive mesh refinement, multi-physics coupling with independent meshes/solvers, arbitrary Lagrangian--Eulerian methods, and overset meshes, where interpolation can introduce violations of physical constraints and, ultimately, the failure of the solver. In these situations, it is necessary to be able to bound the continuous extrema of these polynomial approximations, such as to ensure that interpolated quantities always retain the necessary bounds-preserving properties. 

For both of these applications, the typical approach is to use bounded basis functions in the finite element representation, namely Bernstein polynomials. The positivity and boundedness of these polynomials results in the convex hull property of the basis, where the extrema of any polynomial are bounded by the minimum and maximum coefficients of its Bernstein basis representation~\citep{leroy2008certificates}. This yields a straightforward approach for mesh validity checks and bounds-preserving limiters, where one can simply evaluate the minimum (or maximum) coefficient of the relevant Bernstein representation to yield a lower bound on the true minimum (or maximum) of the polynomial. This approach is widely used in the literature to validate mesh elements (see, for example, \citep{luo2002p, coppeans2024anisotropic, rochery2025metris, johnen2013geometrical, Johnen2017, Mandad2020, Pietroni2022}) and continuously bounds-preserving limiting (see, for example, \citep{Anderson2017, Glaubitz2019,Hajduk2021}). Alternate approaches based on techniques such as sum-of-squares relaxation~\citep{Lasserre2007, marschner2020hexahedral, Despres2020}, other bounded bases~\citep{Hussain2008,Butt1993,DeRossi2017}, and nonlinear optimization~\citep{Dzanic2024, Dzanic2024b,Chen2025} also exist, although typically with higher algorithm complexity and sometimes without strict guarantees on boundedness.

However, the standard methods based on Bernstein representations suffer from two distinct drawbacks, the first being that the transformation from various polynomial representations (e.g., nodal interpolatory bases) to Bernstein polynomials can be ill-conditioned and cause numerical difficulties. More importantly, the bounds one yields from the evaluating the Bernstein coefficients are often ``loose'' in the sense that the bound and the true extrema can differ drastically, and it is very common, for example, to have negative Bernstein basis coefficients even when a polynomial is strictly positive. As such, it is often necessary to subdivide/refine the elements in question a number of times to yield an acceptable level of accuracy in the bounds. This issue is further compounded by the fact that the Bernstein basis is not nodal and its coefficients do not yield any information about the locality of any extrema (i.e., it is not possible to narrow down the location of the extrema within an element without further numerical effort). Consequently, while Bernstein-based approaches provide a systematic means of obtaining bounds, their practical utility is often limited by the conservativeness of these bounds and the computational difficulties associated with obtaining them.

In this work, we introduce a novel approach to computing bounds of high-order polynomials stemming from finite element approximations. 
Broadly inspired by the technique introduced in \citet{Mittal2025} for general field evaluation in high-order finite element methods,  we propose a constrained optimization approach to precompute piecewise linear bounding boxes for the finite element basis functions. These bounding boxes can then be linearly combined to locally bound the extrema of any polynomial that can be recovered from a combination of these basis functions, the accuracy of which can be furthered improved by a simple basis transformation. The proposed approach can be applied to any basis functions/elements, can bound extrema to any desired level of accuracy, and can be efficiently evaluated on-the-fly. We show the applicability of this approach in mesh validity checks and mesh optimization for high-order curved meshes and continuously bounds-preserving limiting for high-order finite element approximations of hyperbolic systems.

The remainder of this manuscript is organized as follows. In \cref{sec:methodology}, we introduce the technique for forming bounding boxes, the optimization process for computing optimal bounding boxes for the basis functions, and an overview of the applications of such approaches to mesh validity checks and bounds-preserving limiting. The results of numerical experiments for these applications are then shown in \cref{sec:results}, followed by conclusions in \cref{sec:conclusion}. 
\section{Methodology}\label{sec:methodology}
Consider a polynomial approximation on the closed subdomain $\Omega$ of the form 
\begin{equation}
    u_h(\mathbf{x}) = \sum_{i=1}^{N} u_i\phi_i(\mathbf{x})\subset V_h, \quad \mathbf{x} \in \Omega,
\end{equation}
where $\phi_i(\mathbf{x})$ are a set of $N$ basis functions of maximal order $p$, $u_i$ are their associated basis coefficients, and $V_h$ is the polynomial space spanned by the basis functions. We use the notation $\mathbb P_p$ to denote a polynomial space of maximal order $p$. This form is commonplace in finite element approximations of partial differential equations, where the subdomain $\Omega$ is an arbitrary element within an arbitrary mesh and $u_h(\mathbf{x})$ approximates some solution within that element. In the one-dimensional case, we take $\Omega = [-1,1]$ and let $\phi_i(\mathbf{x})$ represent finite element basis functions of order $p = N-1$ with $u_i$ as the corresponding degrees of freedom. We do not impose any constraints on the type of basis functions used (i.e., they can be nodal or modal). The goal of this work is to find an efficient technique for bounding the extrema of the high-order polynomial $u_h(\mathbf{x})$ in the form of
\begin{equation}
    u_{\min} \leq u_h(\mathbf{x}) \leq u_{\max} \ \forall \ \mathbf{x} \in \Omega,
\end{equation}
with the accuracy of the bounding method dictated by the ``tightness'' of the bounds with respect to the true extrema, computed as
\begin{equation}
    \left | u_{\min} - \underset{\mathbf{x} \in \Omega}{\min}\ u(\mathbf{x}) \right | \quad\quad \mathrm{and} \quad\quad \left | u_{\max} - \underset{\mathbf{x} \in \Omega}{\max}\ u(\mathbf{x}) \right |.
\end{equation}
For brevity, we neglect the notation $\forall \ \mathbf{x} \in \Omega$ from this point onwards, but it should be understood that the formulations to be presented all operate within some arbitrary element $\Omega$.

\begin{figure}[htbp!]
        \centering
        \adjustbox{width=0.48\linewidth, valign=b}{\begin{tikzpicture}[spy using outlines={rectangle, height=3cm,width=2.5cm, magnification=3, connect spies}]
    \begin{axis}
    [
        axis line style={latex-latex},
        axis y line=left,
        axis x line=left,
        clip mode=individual,
        xlabel = {$x$},
        ylabel = {$u$},
        xmin = -1, xmax = 1,
        ymin = -0.5, ymax = 0.8,
        legend cell align={left},
        legend style={at={(0.97, 0.03)}, anchor=south east},
        x tick label style={/pgf/number format/.cd, fixed, fixed zerofill, precision=1, /tikz/.cd},
        y tick label style={/pgf/number format/.cd, fixed, precision=2, /tikz/.cd},
    ]

        \addplot[color=black, style={very thick},  mark=none, 
        domain=-1:1, samples=100, name path=A]{(-0.5*x^5 + 1*x^4 + 0.1*x^3 - 0.5*x^2 + 0.25*x )};  
        \addlegendentry{$u_h(x)$};
        
        \node [right]  at (axis cs:-1.01,-.2) {\textcolor{PlotColor2}{\footnotesize$(\eta_0, q_0^-)$}};
        \node [right]  at (axis cs:-0.65,-0.35) {\textcolor{PlotColor2}{\footnotesize$(\eta_1, q_1^-)$}};
        \node [right]  at (axis cs:-0.13,-.12) {\textcolor{PlotColor2}{\footnotesize$(\eta_2, q_2^-)$}};
        \node [right]  at (axis cs:0.37,-.02) {\textcolor{PlotColor2}{\footnotesize$(\eta_3, q_3^-)$}};
        \node [right]  at (axis cs:0.9,0.18) {\textcolor{PlotColor2}{\footnotesize$(\eta_4, q_4^-)$}};
        
        \addplot[color=PlotColor2, style={semithick}, mark=*, mark options={scale=0.5}, name path=B] coordinates {
        (-1.0,0.2)
        (-0.6,-0.3)
        (   0,-0.05)
        ( 0.6,0.05)
        ( 1.0,0.25) };
        \addlegendentry{$L^{\boldsymbol{\eta}, \mathbf{q}}(x)$};

        \node [right]  at (axis cs:-1.01,0.72) {\textcolor{PlotColor1}{\footnotesize$(\eta_0, q_0^+)$}};
        \node [right]  at (axis cs:-0.65,-0.05) {\textcolor{PlotColor1}{\footnotesize$(\eta_1, q_1^+)$}};
        \node [right]  at (axis cs:-0.13,0.08) {\textcolor{PlotColor1}{\footnotesize$(\eta_2, q_2^+)$}};
        \node [right]  at (axis cs:0.37,0.16) {\textcolor{PlotColor1}{\footnotesize$(\eta_3, q_3^+)$}};
        \node [right]  at (axis cs:0.9,0.46) {\textcolor{PlotColor1}{\footnotesize$(\eta_4, q_4^+)$}};

        \addplot[color=PlotColor1, style={semithick}, mark=*, mark options={scale=0.5}, name path=B] coordinates {
        (-1.0,0.7)
        (-0.6,-0.15)
        (   0,0.02)
        ( 0.6,0.1)
        ( 1.0,0.4) };
        \addlegendentry{$U^{\boldsymbol{\eta}, \mathbf{q}}(x)$};
        
        \draw[->, style={semithick, dashed}, color=PlotColor2](-0.86, -.17)--(-0.98, 0.15);
    
    \end{axis}
\end{tikzpicture}}
        \caption{Example of $C^0$ upper/lower linear bounding functions for an arbitrary polynomial $u_h(x)$ defined by control nodes $\boldsymbol{\eta}$ and control values $\mathbf{q}$.}
        \label{fig:bound_example}  
    \end{figure}
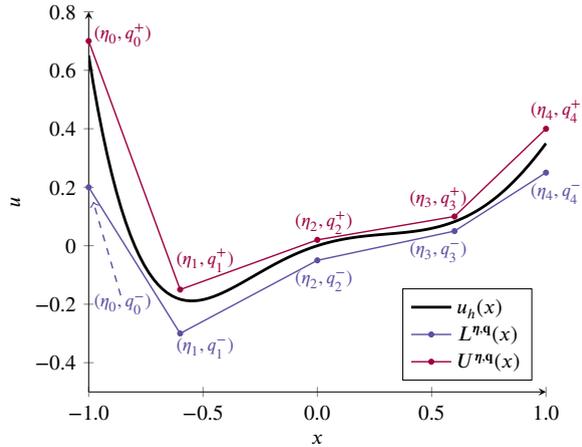

The overarching techniques in this work broadly rely on modifications of the approach introduced in \citet{Mittal2025} for general field evaluation in high-order finite element methods. We present here a brief overview of the approach. Let $\boldsymbol{\eta} \in \Omega$ be some set of $M$ \emph{control nodes} and $\mathbf{q}$ be the associated \emph{control values} of these nodes. For the moment, we assume that the nodes $\boldsymbol{\eta}$ are fixed arbitrarily and consider their optimal positioning later. We define $L^{\boldsymbol{\eta}, \mathbf{q}}(\mathbf{x})$ and $U^{\boldsymbol{\eta}, \mathbf{q}}(\mathbf{x})$ as $C^0$ interpolation functions that linearly interpolate between the control node/value pairs in $\{ \boldsymbol{\eta}, \mathbf{q}\}$. It can be easily seen that if sets of control values $\mathbf{q}_-$ and $\mathbf{q}_+$ are found such that
\begin{equation}\label{eq:condition}
    L^{\boldsymbol{\eta}, \mathbf{q}_-}(\mathbf{x}) \leq u_h(\mathbf{x}) \leq U^{\boldsymbol{\eta}, \mathbf{q}_+}(\mathbf{x}),
\end{equation}
then $u_{\max} = \max \mathbf{q}_+$ and $u_{\min} = \min \mathbf{q}_-$ are upper and lower bounds on $u_h(\mathbf{x})$, respectively. An example visualization of this is shown in \cref{fig:bound_example}. However, finding a set of control nodes/values which guarantees that \cref{eq:condition} is satisfied is non-trivial for arbitrary high-order polynomials and generally requires the solution of a nonlinear optimization problem. For practical applications, it is not feasible to compute these values on-the-fly for arbitrary polynomials, and a more sophisticated approach is required. 

The proposed method in this work relies on two main components. The first component exploits the linearity of the polynomial $u_h(\mathbf{x})$ with respect to the basis functions $\phi(\mathbf{x})$ \emph{which are fixed}. Therefore, for each basis function $\phi_i(\mathbf{x})$, if one precomputes the respective upper/lower control values $\mathbf{q}^i_+$ and $\mathbf{q}^i_-$ such that
\begin{equation}
    L^i(\mathbf{x}) \leq \phi_i(\mathbf{x}) \leq U^i(\mathbf{x}),
\end{equation}
where we use the shorthand notation $L^i(\mathbf{x}) = L^{\boldsymbol{\eta}, \mathbf{q}^i_-}(\mathbf{x})$ and $U^i(\mathbf{x}) = U^{\boldsymbol{\eta}, \mathbf{q}^i_+}(\mathbf{x})$, then one can directly compute linear bounding functions for an arbitrary $u_h(\mathbf{x})$ as 
\begin{subequations}\label{eq:bbox}
    \begin{align}
    L(\mathbf{x}) &= \sum_{i=1}^{N} \min \left (  u_i L^i(\mathbf{x}),\ u_i U^i(\mathbf{x}) \right )\leq u_h(\mathbf{x}),\\
    U(\mathbf{x}) &= \sum_{i=1}^{N} \max \left (  u_i L^i(\mathbf{x}),\ u_i U^i(\mathbf{x}) \right )\geq u_h(\mathbf{x})  .
\end{align}
\end{subequations}
Here, we use the property that
\begin{subequations}\label{eq:scale}
    \begin{align}
     u_i L^i(\mathbf{x}) \leq u_i \phi_i(\mathbf{x}) \leq u_i U^i(\mathbf{x}), \quad \text{if } &u_i > 0\\
     u_i U^i(\mathbf{x}) \leq u_i \phi_i(\mathbf{x}) \leq u_i L^i(\mathbf{x}), \quad \text{if } &u_i < 0.
    \end{align}
\end{subequations}
The control nodes values for these bounding functions, which can be simply computed by summing the respective minima/maxima of the control nodes values for each basis function multiplied by their basis coefficients, yield an algorithmically straightforward approach to computing bounds for an arbitrary polynomial $u_h(\mathbf{x})$ at each control node $\eta_j$ as
\begin{subequations}\label{eq:1D_bounds}
    \begin{align}
    L(\eta_j) &= \sum_{j=1}^N \min(u_i q_{ij}^-,u_i q_{ij}^+),\\
    U(\eta_j) &= \sum_{j=1}^N \max(u_i q_{ij}^-,u_i q_{ij}^+).
    \end{align}
\end{subequations}
We use the notation here that $q_{ij}$ is the value at the control node $\eta_j$ for the basis function $\phi_i$.

However, a naive implementation of the above method results in bounds that are relatively loose. The second component of the proposed method relies on a transformation of the polynomial $u_h(\mathbf{x})$, which helps to tighten these bounds. As the tightness of the bounds is related to the magnitude of the coefficients $u_i$, we consider a transformation of the form
\begin{equation}
    u_h(\mathbf{x}) = \sum_{i=1}^{N} u_i\phi_i(\mathbf{x}) = u_{LO}(\mathbf{x})+ \sum_{i=1}^{N} u_i'\phi_i(\mathbf{x}),
\end{equation}
where we separate $u_h(\mathbf{x})$ into a low-order, \emph{linear} portion ($u_{LO}(\mathbf{x})$) and a high-order portion ($u_h'(\mathbf{x})$), the latter for which we use the shorthand notation $u_h'(\mathbf{x}) = \sum_{i=1}^{N} u_i'\phi_i(\mathbf{x})$. The idea here is to choose the linear portion such that the high-order fluctuations $u_h'(\mathbf{x})$ are relatively ``small''. The method then attempts to bound the high-order fluctuations $u_h'(\mathbf{x})$ instead and superimpose them on the bounds of $u_{LO}(\mathbf{x})$, which can be exactly computed trivially as
\begin{subequations}\label{eq:1D_bounds_offset}
    \begin{align}
    L(\eta_j) = u_{LO}(\eta_j) + L'(\eta_j) =  u_{LO}(\eta_j) +  \sum_{i=1}^N \min \left (u^{'}_i q_{ij}^-,u^{'}_i q_{ij}^+ \right),\\
    U(\eta_j) = u_{LO}(\eta_j) + U'(\eta_j) =  u_{LO}(\eta_j) +\sum_{i=1}^N \max \left (u^{'}_i q_{ij}^-,u^{'}_i q_{ij}^+ \right).
    \end{align}
\end{subequations}

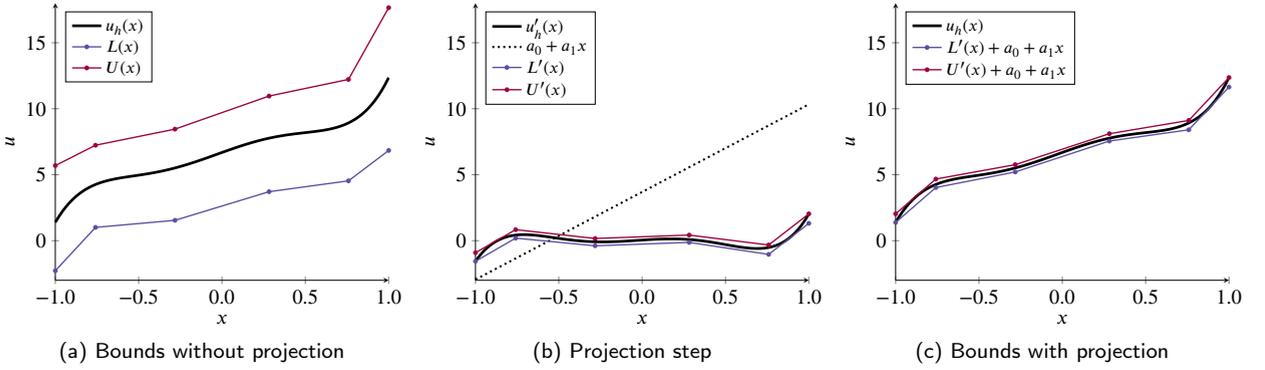
\begin{figure}[htbp!]
        \centering
        \subfloat[Bounds without projection]{
        \adjustbox{width=0.33\linewidth, valign=b}{\begin{tikzpicture}[spy using outlines={rectangle, height=3cm,width=2.5cm, magnification=3, connect spies},
font=\large]
    \begin{axis}
    [
        axis line style={latex-latex},
        axis y line=left,
        axis x line=left,
        clip mode=individual,
        xlabel = {$x$},
        ylabel = {$u$},
        xmin = -1, xmax = 1,
        ymin = -3, ymax = 18,
        legend cell align={left},
        legend style={font=\normalsize, at={(0.03, 0.97)}, anchor=north west},
        x tick label style={/pgf/number format/.cd, fixed, fixed zerofill, precision=1, /tikz/.cd},
        y tick label style={/pgf/number format/.cd, fixed, precision=2, /tikz/.cd},
    ]
        
        \addplot[color=black, style={ultra thick},  mark=none, 
        domain=-1:1, samples=100, name path=A]{(8.2*x^5 + 0.8283*x^4 - 7.234*x^3 - 0.6395*x^2 + 4.52*x + 6.685)}; 
        \addlegendentry{$u_h(x)$};
        
        \addplot[color=PlotColor2, style={thick}, mark=*, mark options={scale=0.5}, name path=B] coordinates {
        (-1.0, -2.2772375208777533)
        (-0.758910868794397, 1.0072495125532994)
        (-0.282320699724306, 1.5428183784258445)
        (0.282320699724306, 3.712039617093923)
        (0.758910868794397, 4.534960789814169)
        (1.0, 6.835982433188313) };
        \addlegendentry{$L(x)$};
        
        \addplot[color=PlotColor1, style={thick}, mark=*, mark options={scale=0.5}, name path=B] coordinates {
        (-1.0, 5.691299775354847)
        (-0.758910868794397, 7.233404907603101)
        (-0.282320699724306, 8.453406789166008)
        (0.282320699724306, 10.961869632817894)
        (0.758910868794397, 12.221986122226056)
        (1.0, 17.658928457457765)
         };
        \addlegendentry{$U(x)$};

    \end{axis}
\end{tikzpicture}}}
        \subfloat[Projection step]{
        \adjustbox{width=0.33\linewidth, valign=b}{\begin{tikzpicture}[spy using outlines={rectangle, height=3cm,width=2.5cm, magnification=3, connect spies},
font=\large]
    \begin{axis}
    [
        axis line style={latex-latex},
        axis y line=left,
        axis x line=left,
        clip mode=individual,
        xlabel = {$x$},
        ylabel = {$u$},
        xmin = -1, xmax = 1,
        ymin = -3, ymax = 18,
        legend cell align={left},
        legend style={font=\normalsize, at={(0.03, 0.97)}, anchor=north west},
        x tick label style={/pgf/number format/.cd, fixed, fixed zerofill, precision=1, /tikz/.cd},
        y tick label style={/pgf/number format/.cd, fixed, precision=2, /tikz/.cd},
    ]
        
        \addplot[color=black, style={ultra thick},  mark=none, 
        domain=-1:1, samples=100, name path=A]{(8.2*x^5 + 0.8283*x^4 - 7.234*x^3 - 0.6395*x^2 + 0.8262*x + 0.04752)}; 
        \addlegendentry{$u_h'(x)$};
        
        \addplot[color=black, style={very thick, dotted},  mark=none, 
        domain=-1:1, samples=100, name path=A]{(6.63726*x + 3.69344)}; 
        \addlegendentry{$a_0 + a_1x$};
        
        \addplot[color=PlotColor2, style={thick}, mark=*, mark options={scale=0.5}, name path=B] coordinates {
        (-1.0, -1.5607611846135359)
(-0.758910868794397, 0.19248728603487836)
(-0.282320699724306, -0.3864933010750274)
(0.282320699724306, -0.12899389527173066)
(0.758910868794397, -1.0387531745716694)
(1.0, 1.3123484580633376) };
        \addlegendentry{$L'(x)$};
        
        \addplot[color=PlotColor1, style={thick}, mark=*, mark options={scale=0.5}, name path=B] coordinates {
        (-1.0, -0.9106617846132419)
(-0.758910868794397, 0.8348116522083724)
(-0.282320699724306, 0.1662622564730714)
(0.282320699724306, 0.425684356750415)
(0.758910868794397, -0.3268597416264241)
(1.0, 2.0344895165203165) };
        \addlegendentry{$U'(x)$};

    \end{axis}
\end{tikzpicture}}}
        \subfloat[Bounds with projection]{
        \adjustbox{width=0.33\linewidth, valign=b}{\begin{tikzpicture}[spy using outlines={rectangle, height=3cm,width=2.5cm, magnification=3, connect spies},
font=\large]
    \begin{axis}
    [
        axis line style={latex-latex},
        axis y line=left,
        axis x line=left,
        clip mode=individual,
        xlabel = {$x$},
        ylabel = {$u$},
        xmin = -1, xmax = 1,
        ymin = -3, ymax = 18,
        legend cell align={left},
        legend style={font=\normalsize, at={(0.03, 0.97)}, anchor=north west},
        x tick label style={/pgf/number format/.cd, fixed, fixed zerofill, precision=1, /tikz/.cd},
        y tick label style={/pgf/number format/.cd, fixed, precision=2, /tikz/.cd},
    ]
        
        \addplot[color=black, style={ultra thick},  mark=none, 
        domain=-1:1, samples=100, name path=A]{(8.2*x^5 + 0.8283*x^4 - 7.234*x^3 - 0.6395*x^2 + 4.52*x + 6.685)}; 
        \addlegendentry{$u_h(x)$};
        
        \addplot[color=PlotColor2, style={thick}, mark=*, mark options={scale=0.5}, name path=B] coordinates {
        (-1.0, 1.3830512983152548)
        (-0.758910868794397, 4.026750397594958)
        (-0.282320699724306, 5.208031724930456)
        (0.282320699724306, 7.551005853613262)
        (0.758910868794397, 8.40150848875873)
        (1.0, 11.643060750025024) };
        \addlegendentry{$L'(x) + a_0 + a_1 x$};

        \addplot[color=PlotColor1, style={thick}, mark=*, mark options={scale=0.5}, name path=B] coordinates {
        (-1.0, 2.033150698315549)
    (-0.758910868794397, 4.669074763768451)
    (-0.282320699724306, 5.760787282478556)
    (0.282320699724306, 8.105684105635408)
    (0.758910868794397, 9.113401921703973)
    (1.0, 12.365201808482002)
         };
        \addlegendentry{$U'(x) + a_0 + a_1 x$};

    \end{axis}
\end{tikzpicture}}}
        \caption{Example of the bounding approach directly applied to a high-order polynomial (left), a visualization of the projection step (middle), and the bounding approach with the projection step (right). }
        \label{fig:proj_example}
    \end{figure}

Per \citet{Mittal2025} and \cite{gslib-github}, in the one-dimensional case, the linear portion can be represented as $u_{LO}(x) = a_0 + a_1 x$ (i.e., the equivalent $\mathbb P_1$ basis for the given element type). The coefficients $a_0$ and $a_1$ are calculated from the $L^2$ projection of the polynomial $u_h(\mathbf{x})$ onto the $\mathbb P_1$ subspace, which correspond to the zeroth/first-order modal coefficients for an orthogonal basis computed with respect to the unit measure (e.g., Legendre basis coefficients for tensor-product elements):
\begin{subequations}\label{eq:projcoeff}
    \begin{align}
    a_0 &= \frac{1}{2}\int_\Omega u_h(x)\ \mathrm{d}x,\\
    a_1 &= \frac{3}{2}\int_\Omega x u_h(x)\ \mathrm{d}x.
    \end{align}
\end{subequations}
A visualization of how the $L^2$ projection to $\mathbb P_1$ compacts the bounds is presented in \cref{fig:proj_example}. The initially loose bounds on a high-order polynomial ($L(x) \leq u_h(x) \leq U(x)$) are converted into a linear component ($u_{LO}(x) = a_0 + a_1 x$) and tight bounds on the high-order fluctuations of the polynomial ($L'(x) \leq u_h'(x) \leq U'(x)$). This modification also gives the benefit of ensuring that the bounding approach is invariant to shifts/linear scalings of the polynomial. Extensions to higher dimensions follow the same approach, although for elements with tensor-product structures, a more efficient method can be obtained by decomposing the problem into a series of one-dimensional problems (see \cref{app:higherdim}).

\subsection{Calculating bounding boxes}\label{ssec:bbox}
The above approach requires calculating a set of control nodes/value pairs $\{ \boldsymbol{\eta}, \mathbf{q}\}$ to create a bounding box for each basis function. In \citet{Mittal2025}, this was approximated by using the values and gradients of the $N$ basis functions (in this case, the nodal interpolating basis functions at the Gauss--Lobatto nodes) at $M$ Chebyshev nodes. This approach, however, is heuristic-based, and the bounding property in \cref{eq:condition} is only shown empirically for certain values of $M > N$, with the minimum necessary resolution dependent on the choice of control nodes. One of the novelties of this work is to present an approach to precompute the bounding functions $L^i(\mathbf{x})$ and $U^i(\mathbf{x})$ which: i) can be applied to any basis function $\phi_i(\mathbf{x})$; ii) can be applied to any set of control nodes $\phi_i(\mathbf{x})$; and iii) are guaranteed to bound the basis functions (i.e., satisfy \cref{eq:condition}) for any number of control nodes $M \geq 2$. This results in a provably robust approach to bound extrema of arbitrary polynomial approximations in finite element methods. We present this first in terms of an arbitrary set of control nodes $\boldsymbol{\eta}$ and later discuss how these control nodes can be selected. 

The task of finding a linear bounding function that ``tightly'' bounds a basis function can be framed as a constrained optimization problem. For a given set of control nodes $\boldsymbol{\eta}$ and basis function $\phi_i(\mathbf{x})$, we seek to find the solution of the optimization problem 
\begin{equation}
    \mathbf{q}^* = \underset{\mathbf{q}}{\arg \min}\ f(\mathbf{q}) \quad s.t.\quad  g(\mathbf{q}, \mathbf{x}) \geq 0 ,
\end{equation}
where
\begin{subequations}\label{eq:obj1}
    \begin{align}
    f(\mathbf{q}) &= \left \| U^i_{\boldsymbol{\eta}, \mathbf{q}}(\mathbf{x}) - \phi_i(\mathbf{x})\right \|_{2, \Omega},\\
    g(\mathbf{q}, \mathbf{x}) &= U^i_{\boldsymbol{\eta}, \mathbf{q}}(\mathbf{x}) - \phi_i(\mathbf{x}).
    \end{align}
\end{subequations}
Note that the constraint functional is framed in terms of an upper bounding function, but an identical formulation can be made for the lower bounding function by negating the inequality of the constraint functional. One can even take this a step further and try to find an \emph{optimal} set of $M$ control nodes $\boldsymbol{\eta}$ over a set of $N$ given basis functions as 
\begin{equation}
    \{\boldsymbol{\eta}^*, \mathbf{q}^*\} = \underset{\{\boldsymbol{\eta}, \mathbf{q}\}}{\arg \min}\ f(\boldsymbol{\eta}, \mathbf{q}) \quad s.t. \quad g^i(\boldsymbol{\eta}, \mathbf{q}, \mathbf{x}) \geq 0, 
\end{equation}
for all $i$ in $\{1, \ldots, M \}$, where
\begin{subequations}\label{eq:obj2}
    \begin{align}
    f(\boldsymbol{\eta}, \mathbf{q}) &= \sum_{i=1}^M \left \|U^i_{\boldsymbol{\eta}, \mathbf{q}}(\mathbf{x}) - \phi_i(\mathbf{x})\right \|_{2, \Omega},\\
    g^i(\mathbf{q}, \mathbf{x}) &= U^i_{\boldsymbol{\eta}, \mathbf{q}}(\mathbf{x}) - \phi_i(\mathbf{x}).
    \end{align}
\end{subequations}
We refer to the solution of the former as the optimal bounds for a given set of points and the latter as the optimal bounds on the optimal points. The optimization process and implementation is further described in \cref{ssec:opt}.

\begin{figure}[htbp!]
        \centering
        \hspace{0.25em}
        \adjustbox{width=0.32\linewidth, valign=b}{\begin{tikzpicture}[spy using outlines={rectangle, height=3cm,width=2.5cm, magnification=3, connect spies},
font=\large]
    \begin{axis}
    [
        axis line style={latex-latex},
        axis y line=left,
        axis x line=left,
        clip mode=individual,
        xlabel = {$x$},
        ylabel = {$\phi_1(x)$},
        xmin = -1,  xmax = 1,
        legend cell align={left},
        legend style={font=\normalsize, at={(0.97, 0.97)}, anchor=north east},
        x tick label style={/pgf/number format/.cd, fixed, fixed zerofill, precision=1, /tikz/.cd},
        y tick label style={/pgf/number format/.cd, fixed, precision=2, /tikz/.cd},	
        major grid style={line width=.2pt,draw=gray!50},
        minor x tick num=4,
        minor y tick num=4,
    ]
        
        \addplot[color=black, style={ultra thick}] table[x=x, y=phi1, col sep=comma]{./figs/data/example_lobatto_N4_M5.csv};
        \addlegendentry{$\phi(x)$};
        
        \addplot[color=PlotColor2, style={very thick}] table[x=x, y=l1, col sep=comma]{./figs/data/example_lobatto_N4_M5.csv};
        \addlegendentry{$L(x)$};
        
        \addplot[color=PlotColor1, style={very thick}] table[x=x, y=u1, col sep=comma]{./figs/data/example_lobatto_N4_M5.csv};
        \addlegendentry{$U(x)$};

    \end{axis}
\end{tikzpicture}}
        \adjustbox{width=0.32\linewidth, valign=b}{\begin{tikzpicture}[spy using outlines={rectangle, height=3cm,width=2.5cm, magnification=3, connect spies},
font=\large]
    \begin{axis}
    [
        axis line style={latex-latex},
        axis y line=left,
        axis x line=left,
        clip mode=individual,
        xlabel = {$x$},
        ylabel = {$\phi_1(x)$},
        xmin = -1,  xmax = 1,
        legend cell align={left},
        legend style={font=\scriptsize, at={(0.97, 0.97)}, anchor=north east},
        x tick label style={/pgf/number format/.cd, fixed, fixed zerofill, precision=1, /tikz/.cd},
        y tick label style={/pgf/number format/.cd, fixed, precision=2, /tikz/.cd},	
        major grid style={line width=.2pt,draw=gray!50},
        minor x tick num=4,
        minor y tick num=4,
    ]
        
        \addplot[color=black, style={ultra thick}] table[x=x, y=phi1, col sep=comma]{./figs/data/example_legendre_N4_M5.csv};
        
        \addplot[color=PlotColor2, style={very thick}] table[x=x, y=l1, col sep=comma]{./figs/data/example_legendre_N4_M5.csv};
        
        \addplot[color=PlotColor1, style={very thick}] table[x=x, y=u1, col sep=comma]{./figs/data/example_legendre_N4_M5.csv};

    \end{axis}
\end{tikzpicture}}
        \adjustbox{width=0.32\linewidth, valign=b}{\begin{tikzpicture}[spy using outlines={rectangle, height=3cm,width=2.5cm, magnification=3, connect spies},
font=\large]
    \begin{axis}
    [
        axis line style={latex-latex},
        axis y line=left,
        axis x line=left,
        clip mode=individual,
        xlabel = {$x$},
        ylabel = {$\phi_1(x)$},
        xmin = -1,  xmax = 1,
        legend cell align={left},
        legend style={font=\scriptsize, at={(0.97, 0.97)}, anchor=north east},
        x tick label style={/pgf/number format/.cd, fixed, fixed zerofill, precision=1, /tikz/.cd},
        y tick label style={/pgf/number format/.cd, fixed, precision=2, /tikz/.cd},	
        major grid style={line width=.2pt,draw=gray!50},
        minor x tick num=4,
        minor y tick num=4,
    ]
        
        \addplot[color=black, style={ultra thick}] table[x=x, y=phi1, col sep=comma]{./figs/data/example_bernstein_N4_M5.csv};
        
        \addplot[color=PlotColor2, style={very thick}] table[x=x, y=l1, col sep=comma]{./figs/data/example_bernstein_N4_M5.csv};
        
        \addplot[color=PlotColor1, style={very thick}] table[x=x, y=u1, col sep=comma]{./figs/data/example_bernstein_N4_M5.csv};

    \end{axis}
\end{tikzpicture}}
        \newline
        \adjustbox{width=0.32\linewidth, valign=b}{\begin{tikzpicture}[spy using outlines={rectangle, height=3cm,width=2.5cm, magnification=3, connect spies},
font=\large]
    \begin{axis}
    [
        axis line style={latex-latex},
        axis y line=left,
        axis x line=left,
        clip mode=individual,
        xlabel = {$x$},
        ylabel = {$\phi_2(x)$},
        xmin = -1,  xmax = 1,
        legend cell align={left},
        legend style={font=\scriptsize, at={(0.97, 0.97)}, anchor=north east},
        x tick label style={/pgf/number format/.cd, fixed, fixed zerofill, precision=1, /tikz/.cd},
        y tick label style={/pgf/number format/.cd, fixed, precision=2, /tikz/.cd},	
        major grid style={line width=.2pt,draw=gray!50},
        minor x tick num=4,
        minor y tick num=4,
    ]
        
        \addplot[color=black, style={ultra thick}] table[x=x, y=phi2, col sep=comma]{./figs/data/example_lobatto_N4_M5.csv};
        
        \addplot[color=PlotColor2, style={very thick}] table[x=x, y=l2, col sep=comma]{./figs/data/example_lobatto_N4_M5.csv};
        
        \addplot[color=PlotColor1, style={very thick}] table[x=x, y=u2, col sep=comma]{./figs/data/example_lobatto_N4_M5.csv};

    \end{axis}
\end{tikzpicture}}
        \adjustbox{width=0.32\linewidth, valign=b}{\begin{tikzpicture}[spy using outlines={rectangle, height=3cm,width=2.5cm, magnification=3, connect spies},
font=\large]
    \begin{axis}
    [
        axis line style={latex-latex},
        axis y line=left,
        axis x line=left,
        clip mode=individual,
        xlabel = {$x$},
        ylabel = {$\phi_2(x)$},
        xmin = -1,  xmax = 1,
        legend cell align={left},
        legend style={font=\scriptsize, at={(0.97, 0.97)}, anchor=north east},
        x tick label style={/pgf/number format/.cd, fixed, fixed zerofill, precision=1, /tikz/.cd},
        y tick label style={/pgf/number format/.cd, fixed, precision=2, /tikz/.cd},	
        major grid style={line width=.2pt,draw=gray!50},
        minor x tick num=4,
        minor y tick num=4,
    ]
        
        \addplot[color=black, style={ultra thick}] table[x=x, y=phi2, col sep=comma]{./figs/data/example_legendre_N4_M5.csv};
        
        \addplot[color=PlotColor2, style={very thick}] table[x=x, y=l2, col sep=comma]{./figs/data/example_legendre_N4_M5.csv};
        
        \addplot[color=PlotColor1, style={very thick}] table[x=x, y=u2, col sep=comma]{./figs/data/example_legendre_N4_M5.csv};

    \end{axis}
\end{tikzpicture}}
        \adjustbox{width=0.32\linewidth, valign=b}{\begin{tikzpicture}[spy using outlines={rectangle, height=3cm,width=2.5cm, magnification=3, connect spies},
font=\large]
    \begin{axis}
    [
        axis line style={latex-latex},
        axis y line=left,
        axis x line=left,
        clip mode=individual,
        xlabel = {$x$},
        ylabel = {$\phi_2(x)$},
        xmin = -1,  xmax = 1,
        legend cell align={left},
        legend style={font=\scriptsize, at={(0.97, 0.97)}, anchor=north east},
        x tick label style={/pgf/number format/.cd, fixed, fixed zerofill, precision=1, /tikz/.cd},
        y tick label style={/pgf/number format/.cd, fixed, precision=2, /tikz/.cd},	
        major grid style={line width=.2pt,draw=gray!50},
        minor x tick num=4,
        minor y tick num=4,
    ]
        
        \addplot[color=black, style={ultra thick}] table[x=x, y=phi2, col sep=comma]{./figs/data/example_bernstein_N4_M5.csv};
        
        \addplot[color=PlotColor2, style={very thick}] table[x=x, y=l2, col sep=comma]{./figs/data/example_bernstein_N4_M5.csv};
        
        \addplot[color=PlotColor1, style={very thick}] table[x=x, y=u2, col sep=comma]{./figs/data/example_bernstein_N4_M5.csv};

    \end{axis}
\end{tikzpicture}}
        \newline
        \adjustbox{width=0.32\linewidth, valign=b}{\begin{tikzpicture}[spy using outlines={rectangle, height=3cm,width=2.5cm, magnification=3, connect spies},
font=\large]
    \begin{axis}
    [
        axis line style={latex-latex},
        axis y line=left,
        axis x line=left,
        clip mode=individual,
        xlabel = {$x$},
        ylabel = {$\phi_3(x)$},
        xmin = -1,  xmax = 1,
        legend cell align={left},
        legend style={font=\scriptsize, at={(0.97, 0.97)}, anchor=north east},
        x tick label style={/pgf/number format/.cd, fixed, fixed zerofill, precision=1, /tikz/.cd},
        y tick label style={/pgf/number format/.cd, fixed, precision=2, /tikz/.cd},	
        major grid style={line width=.2pt,draw=gray!50},
        minor x tick num=4,
        minor y tick num=4,
    ]
        
        \addplot[color=black, style={ultra thick}] table[x=x, y=phi3, col sep=comma]{./figs/data/example_lobatto_N4_M5.csv};
        
        \addplot[color=PlotColor2, style={very thick}] table[x=x, y=l3, col sep=comma]{./figs/data/example_lobatto_N4_M5.csv};
        
        \addplot[color=PlotColor1, style={very thick}] table[x=x, y=u3, col sep=comma]{./figs/data/example_lobatto_N4_M5.csv};

    \end{axis}
\end{tikzpicture}}
        \adjustbox{width=0.32\linewidth, valign=b}{\begin{tikzpicture}[spy using outlines={rectangle, height=3cm,width=2.5cm, magnification=3, connect spies},
font=\large]
    \begin{axis}
    [
        axis line style={latex-latex},
        axis y line=left,
        axis x line=left,
        clip mode=individual,
        xlabel = {$x$},
        ylabel = {$\phi_3(x)$},
        xmin = -1,  xmax = 1,
        legend cell align={left},
        legend style={font=\scriptsize, at={(0.97, 0.97)}, anchor=north east},
        x tick label style={/pgf/number format/.cd, fixed, fixed zerofill, precision=1, /tikz/.cd},
        y tick label style={/pgf/number format/.cd, fixed, precision=2, /tikz/.cd},	
        major grid style={line width=.2pt,draw=gray!50},
        minor x tick num=4,
        minor y tick num=4,
    ]
        
        \addplot[color=black, style={ultra thick}] table[x=x, y=phi3, col sep=comma]{./figs/data/example_legendre_N4_M5.csv};
        
        \addplot[color=PlotColor2, style={very thick}] table[x=x, y=l3, col sep=comma]{./figs/data/example_legendre_N4_M5.csv};
        
        \addplot[color=PlotColor1, style={very thick}] table[x=x, y=u3, col sep=comma]{./figs/data/example_legendre_N4_M5.csv};

    \end{axis}
\end{tikzpicture}}
        \adjustbox{width=0.32\linewidth, valign=b}{\begin{tikzpicture}[spy using outlines={rectangle, height=3cm,width=2.5cm, magnification=3, connect spies},
font=\large]
    \begin{axis}
    [
        axis line style={latex-latex},
        axis y line=left,
        axis x line=left,
        clip mode=individual,
        xlabel = {$x$},
        ylabel = {$\phi_3(x)$},
        xmin = -1,  xmax = 1,
        legend cell align={left},
        legend style={font=\scriptsize, at={(0.97, 0.97)}, anchor=north east},
        x tick label style={/pgf/number format/.cd, fixed, fixed zerofill, precision=1, /tikz/.cd},
        y tick label style={/pgf/number format/.cd, fixed, precision=2, /tikz/.cd},	
        major grid style={line width=.2pt,draw=gray!50},
        minor x tick num=4,
        minor y tick num=4,
    ]
        
        \addplot[color=black, style={ultra thick}] table[x=x, y=phi3, col sep=comma]{./figs/data/example_bernstein_N4_M5.csv};
        
        \addplot[color=PlotColor2, style={very thick}] table[x=x, y=l3, col sep=comma]{./figs/data/example_bernstein_N4_M5.csv};
        
        \addplot[color=PlotColor1, style={very thick}] table[x=x, y=u3, col sep=comma]{./figs/data/example_bernstein_N4_M5.csv};

    \end{axis}
\end{tikzpicture}}
        \newline
        \subfloat[Gauss--Lobatto basis]{
        \adjustbox{width=0.32\linewidth, valign=b}{\begin{tikzpicture}[spy using outlines={rectangle, height=3cm,width=2.5cm, magnification=3, connect spies},
font=\large]
    \begin{axis}
    [
        axis line style={latex-latex},
        axis y line=left,
        axis x line=left,
        clip mode=individual,
        xlabel = {$x$},
        ylabel = {$\phi_4(x)$},
        xmin = -1,  xmax = 1,
        legend cell align={left},
        legend style={font=\scriptsize, at={(0.97, 0.97)}, anchor=north east},
        x tick label style={/pgf/number format/.cd, fixed, fixed zerofill, precision=1, /tikz/.cd},
        y tick label style={/pgf/number format/.cd, fixed, precision=2, /tikz/.cd},	
        major grid style={line width=.2pt,draw=gray!50},
        minor x tick num=4,
        minor y tick num=4,
    ]
        
        \addplot[color=black, style={ultra thick}] table[x=x, y=phi4, col sep=comma]{./figs/data/example_lobatto_N4_M5.csv};
        
        \addplot[color=PlotColor2, style={very thick}] table[x=x, y=l4, col sep=comma]{./figs/data/example_lobatto_N4_M5.csv};
        
        \addplot[color=PlotColor1, style={very thick}] table[x=x, y=u4, col sep=comma]{./figs/data/example_lobatto_N4_M5.csv};

    \end{axis}
\end{tikzpicture}}}
        \subfloat[Gauss--Legendre basis]{
        \adjustbox{width=0.32\linewidth, valign=b}{\begin{tikzpicture}[spy using outlines={rectangle, height=3cm,width=2.5cm, magnification=3, connect spies},
font=\large]
    \begin{axis}
    [
        axis line style={latex-latex},
        axis y line=left,
        axis x line=left,
        clip mode=individual,
        xlabel = {$x$},
        ylabel = {$\phi_4(x)$},
        xmin = -1,  xmax = 1,
        legend cell align={left},
        legend style={font=\scriptsize, at={(0.97, 0.97)}, anchor=north east},
        x tick label style={/pgf/number format/.cd, fixed, fixed zerofill, precision=1, /tikz/.cd},
        y tick label style={/pgf/number format/.cd, fixed, precision=2, /tikz/.cd},	
        major grid style={line width=.2pt,draw=gray!50},
        minor x tick num=4,
        minor y tick num=4,
    ]
        
        \addplot[color=black, style={ultra thick}] table[x=x, y=phi4, col sep=comma]{./figs/data/example_legendre_N4_M5.csv};
        
        \addplot[color=PlotColor2, style={very thick}] table[x=x, y=l4, col sep=comma]{./figs/data/example_legendre_N4_M5.csv};
        
        \addplot[color=PlotColor1, style={very thick}] table[x=x, y=u4, col sep=comma]{./figs/data/example_legendre_N4_M5.csv};

    \end{axis}
\end{tikzpicture}}}
        \subfloat[Bernstein basis]{
        \adjustbox{width=0.32\linewidth, valign=b}{\begin{tikzpicture}[spy using outlines={rectangle, height=3cm,width=2.5cm, magnification=3, connect spies},
font=\large]
    \begin{axis}
    [
        axis line style={latex-latex},
        axis y line=left,
        axis x line=left,
        clip mode=individual,
        xlabel = {$x$},
        ylabel = {$\phi_4(x)$},
        xmin = -1,  xmax = 1,
        legend cell align={left},
        legend style={font=\scriptsize, at={(0.97, 0.97)}, anchor=north east},
        x tick label style={/pgf/number format/.cd, fixed, fixed zerofill, precision=1, /tikz/.cd},
        y tick label style={/pgf/number format/.cd, fixed, precision=2, /tikz/.cd},	
        major grid style={line width=.2pt,draw=gray!50},
        minor x tick num=4,
        minor y tick num=4,
    ]
        
        \addplot[color=black, style={ultra thick}] table[x=x, y=phi4, col sep=comma]{./figs/data/example_bernstein_N4_M5.csv};
        
        \addplot[color=PlotColor2, style={very thick}] table[x=x, y=l4, col sep=comma]{./figs/data/example_bernstein_N4_M5.csv};
        
        \addplot[color=PlotColor1, style={very thick}] table[x=x, y=u4, col sep=comma]{./figs/data/example_bernstein_N4_M5.csv};

    \end{axis}
\end{tikzpicture}}}
        \newline
        \caption{Examples of optimal bounding boxes for $N=4$ Gauss--Lobatto (left), Gauss--Legendre (middle), and Bernstein (right) basis functions with $M=5$ equispaced control nodes. }
        \label{fig:bases}  
    \end{figure}
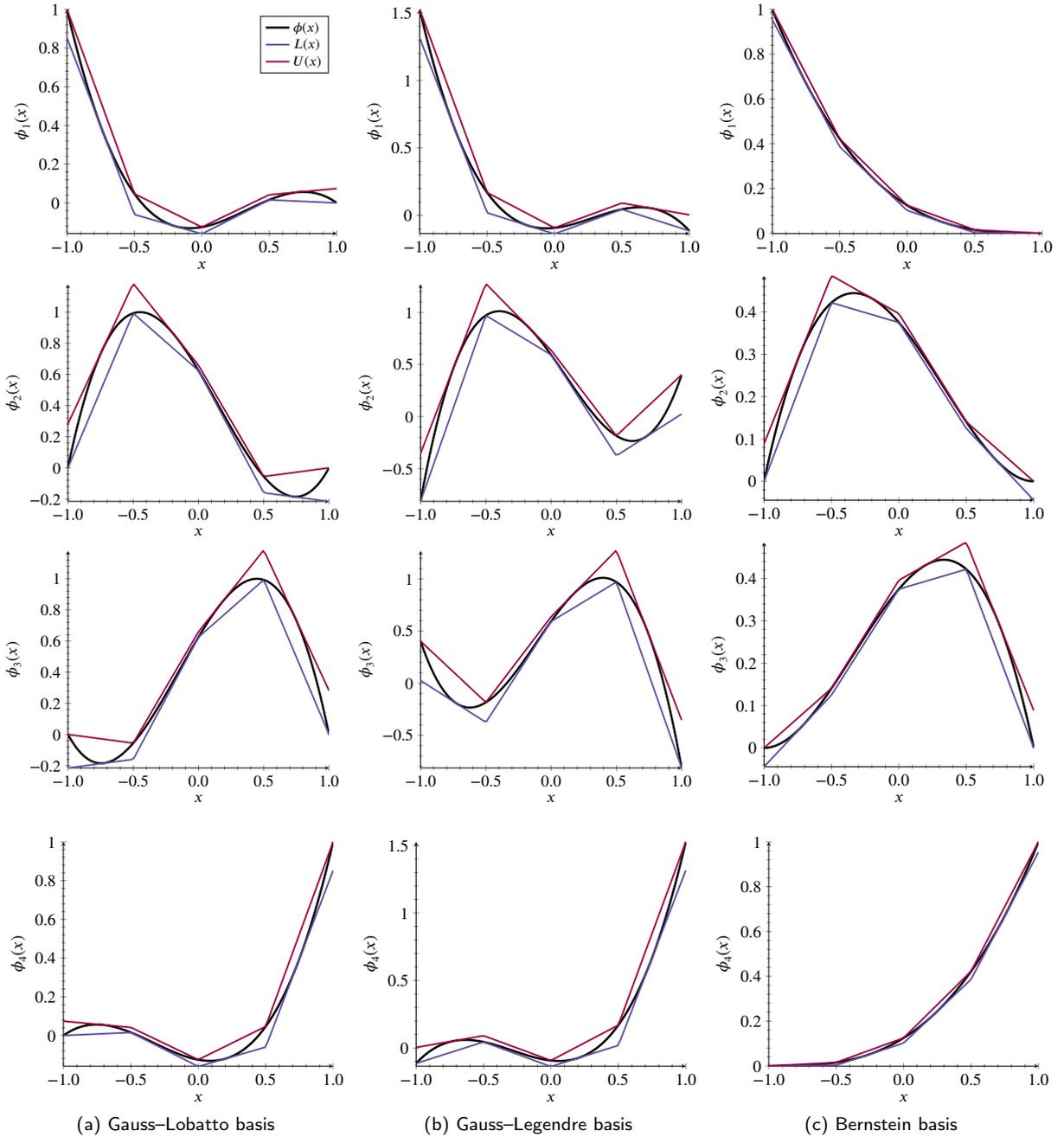

\emph{This proposed approach has the benefit that the bounding properties and optimality of the bounding functions are independent of the choice of basis functions and control points.} Unlike the original approach of \citet{Mittal2025}, the bounding functions are guaranteed to bound any basis function for any value of $M > 1$ (per dimension), and the bounds are optimal in the sense that they minimize the $L^2$ norm of the bounding error. We showcase examples of these computed bounding functions for various bases, including Gauss--Lobatto, Gauss--Legendre, and Bernstein bases, in \cref{fig:bases}. As this optimization process only needs to be performed once given a basis function and a desired number of control nodes $M$, the control nodes/values for a variety of bases can be precomputed and stored for future simulations. Therefore, the actual cost of the bounding technique during simulations for general elements is simply $\mathcal O(NM)$ operations (for computing the summations in \cref{eq:bbox}) for one-dimensional functions. For higher-dimensional functions with tensor-product structures, the problem can be treated as a sequence of one-dimensional problems, which is further explained in \cref{app:higherdim}. For the numerical results in this work, we present $N$ and $M$ in terms of the one-dimensional values for clarity.

\begin{figure}[htbp!]
        \centering
        \subfloat[Base resolution]{
        \adjustbox{width=0.32\linewidth, valign=b}{\includegraphics[]{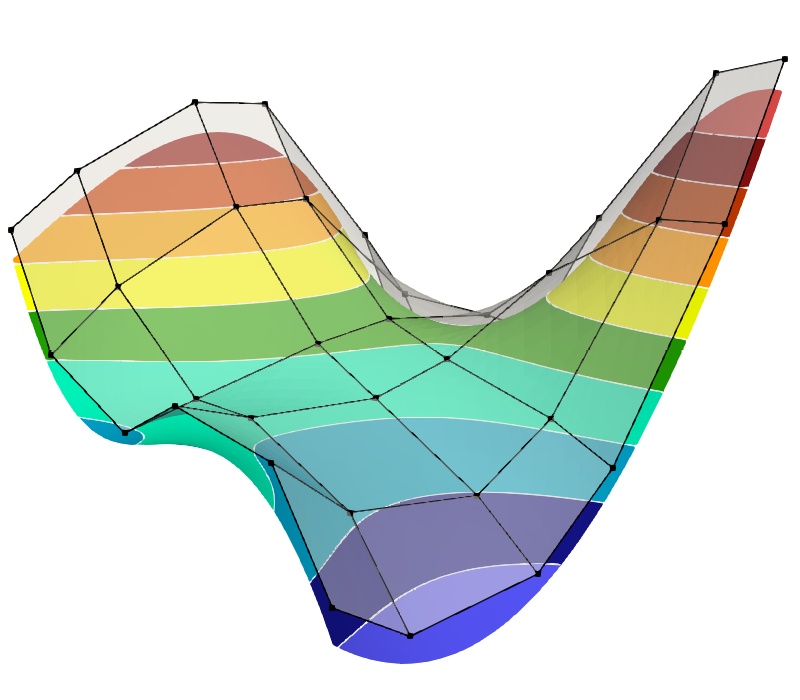}}}
        \hfill
        \subfloat[Increased number of control points]{
        \adjustbox{width=0.32\linewidth, valign=b}{\includegraphics[]{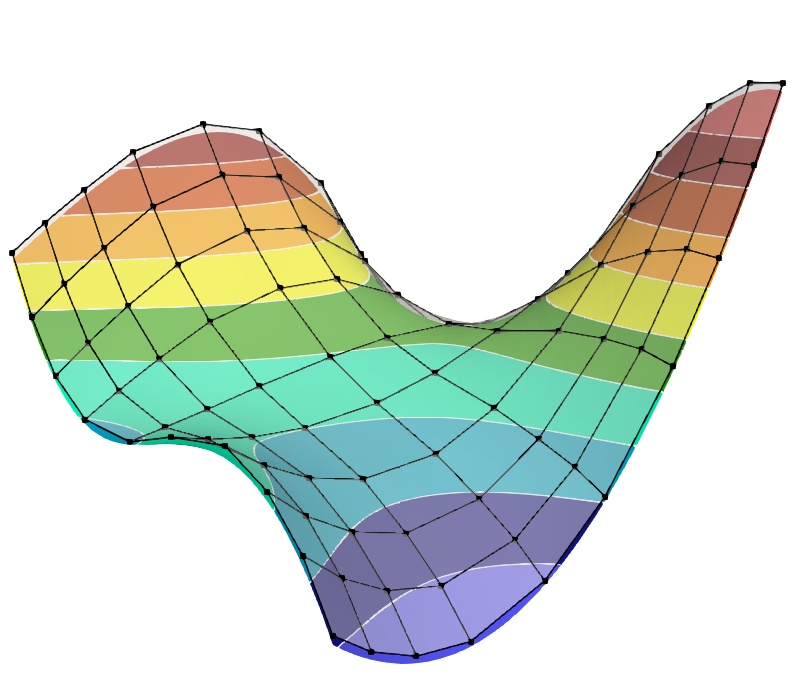}}}
        \hfill
        \subfloat[Increased subdivision levels]{
        \adjustbox{width=0.32\linewidth, valign=b}{\includegraphics[]{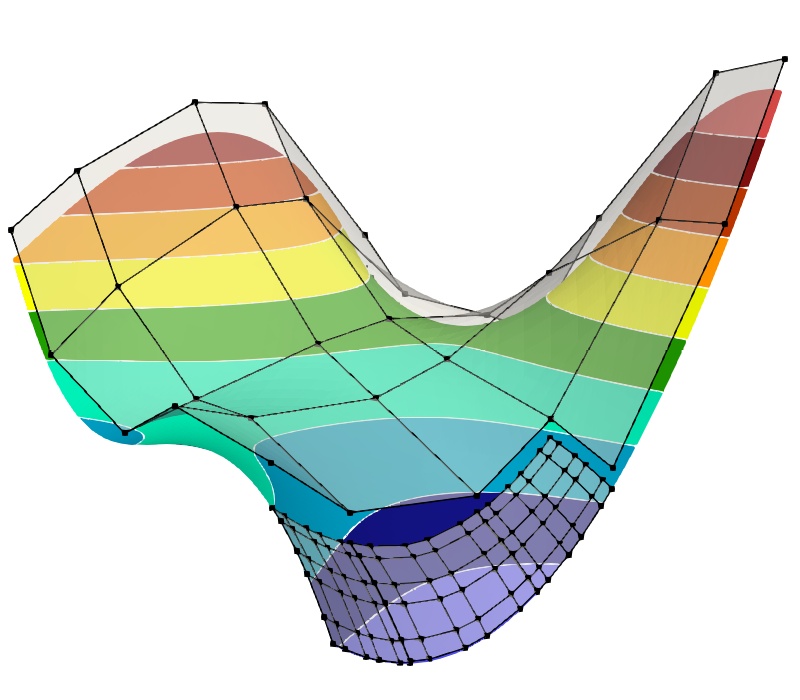}}}
        \caption{Example of an upper bounding surface for a high-order two-dimensional polynomial at some base level of resolution (left), increased resolution through an increased number of control points $M$ (middle), and increased resolution through subdivision (right).}
        \label{fig:2d_schematic}
    \end{figure}

Compared to some approaches that rely on convex hull-type properties for bounding polynomials (i.e., Bernstein polynomial representations), the proposed bounding technique has the benefit of giving localized bounds in terms of the values at the control nodes, which allows for refined estimates of local extrema. In contrast, other bounding approaches in the literature often can only provide bounds for the entire element, which we refer to as global bounds. From this, one can increase the resolution to control the accuracy of the bounds, generally dictated by some desired tolerance level between the local lower and upper bounds. One such approach is to increase the number of control nodes $M$ -- it will later be shown that the bounding error is second-order with respect to $M$. This allows for the user to obtain initially very loose bounds (e.g., using $M = 2$) and then increase resolution as necessary. Alternatively, if one wishes to further refine the bounds in a particular subregion of the element, a common technique is to subdivide the element into sub-elements (similarly to adaptive mesh refinement) and apply the bounding technique to the particular sub-element. Compared to simply increasing $M$, this provides the benefit of completely localizing any refinements in the bounds estimates, although at the expense of more algorithmically intensive operations such as interpolation onto the sub-element and recalculation of the bounds. We note here that compared to global bounds estimates (such as taking the minimum/maximum Bernstein coefficients over the element), the local nature of the proposed bounding approach allows one to essentially skip one level of refinement at the beginning and immediately refine a desired sub-element between control nodes instead of the entire element. We show a visualization of these two refinement approaches in \cref{fig:2d_schematic} and show examples of both in the numerical experiments. 

\begin{figure}[htbp!]
        \centering
        \adjustbox{width=0.32\linewidth, valign=b}{\includegraphics[]{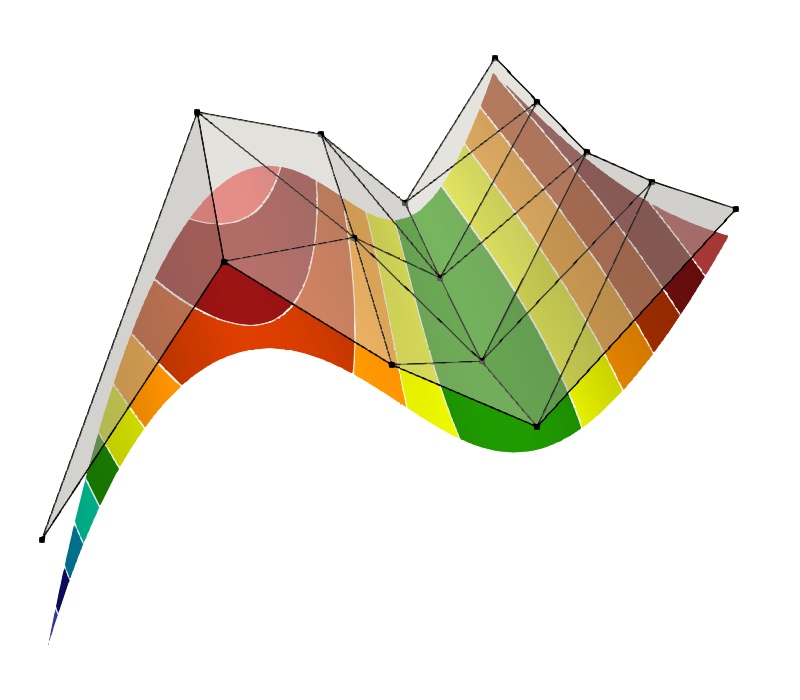}}
        \caption{Example of an upper bounding surface for a high-order two-dimensional polynomial on a triangular element with equispaced control nodes.}
        \label{fig:triangle}
    \end{figure}

While the focus of the numerical experiments in this work is on tensor-product elements, the presented formulation for optimizing the bounding boxes and computing the bounds is applicable to any element type. We showcase an example of the identical formulation applied to bounding a high-order polynomial on a triangle element on a fixed set of equispaced control nodes in \cref{fig:triangle}, where we optimize \cref{eq:obj1} now with respect to the two-dimensional nodal basis functions on the reference triangle. While this optimization is straightforward for a given set of control nodes (albeit at a notably larger computational cost due to the increased dimensionality), the process of computing the optimal control nodes becomes much more difficult. For tensor-product elements, one can simply compute the higher-dimensional bounding boxes as tensor-products of the one-dimensional bounding boxes, and enforcing symmetry on the control node distribution and bounding boxes is straightforward. For non-tensor product elements, the two main complications that arise are: i) one must choose how to tessellate between an arbitrary set of nodes (the optimality of which is an open problem); and ii) constraining the node distributions within the element and enforcing symmetry in the control node distribution is notably more complicated (e.g., symmetry orbits used for quadrature rules can offer some insights~\citep{Witherden2015}). As such, we leave the topic of bounding extrema for simplex elements as a topic of future work.

\subsection{Application to mesh validity checks}
The bounding of polynomial extrema is of interest to problems in meshing, particularly for high-order/curved meshes which deform nonlinearly. In these applications, elements are often represented in terms of a mapping $\mathbf{J}:\boldsymbol{\xi} \to \mathbf{x}$ which transforms the element from a reference space $\boldsymbol{\xi}$ to the physical space $\mathbf{x}$. The entries of this Jacobian transformation matrix $\mathbf{J}$, given as
\begin{equation}
    J_{ij} = \frac{\partial \mathbf{x}_i(\boldsymbol{\xi})}{\partial \boldsymbol{\xi}_j},
\end{equation}
are typically polynomial functions of the reference coordinates $\boldsymbol{\xi}$, and their extrema determine element validity and mesh quality. In particular, the local volume of the element can be inferred from the determinant of the Jacobian $|\mathbf{J}|$, and the validity of the element is conditional on the positivity of this determinant across the entire element (i.e., $|\mathbf{J}| > 0 \ \forall\ \mathbf{x}$ if the element is not inverted). 

For linear simplex elements, this can be verified by simply computing the determinant at the mesh nodes. For higher-order meshes, this cannot be verified in the same manner as the positivity of the determinant at mesh nodes does not imply positivity throughout the entire element, which can cause non-valid elements with regions of negative volume/inversion~\citep{dobrev2019target,sanjaya2020comparison,aparicio2024defining}. Alternatively, one can use Bernstein bases for determinant~\citep{coppeans2024anisotropic,luo2002p,johnen2013geometrical,rochery2025metris}, for which validity can be guaranteed if the minimum basis coefficient is positive, but this approach is suboptimal as it can be computationally expensive and the bounds are relatively loose (i.e., it is quite common to have the minimum basis coefficient be negative even when an element is valid ). As such, one application of the proposed technique is in checking the validity of high-order/curved meshes. In particular, we can represent the determinant of the Jacobian as a polynomial in the reference domain, and it is simple to verify that if a $d$-dimensional mesh is of order $p$ (i.e., the maximal order of $J_{ij}$ is of order $p$), then the maximal order of $|\mathbf{J}|({\xi})$ is $dp-1$ for tensor-product elements and $d(p-1)$ for simplices~\citep{johnen2013geometrical}. Therefore, we can simply treat $|\mathbf{J}|({\xi})$ as a ``solution'' polynomial of that order and apply the proposed bounding technique. 

In particular, if the minimum bound for $|\mathbf{J}|({\xi})$ is positive \emph{at all control nodes}, then we can guarantee that the element is valid. Furthermore, if the maximum bound for $|\mathbf{J}|({\xi})$ is negative \emph{at any control node}, we can guarantee that the element is invalid. If the maximum bound is positive but the minimum bound is negative \emph{at the same control node}, then the element may or may not be valid. In this scenario, we can increase the resolution, either through subdividing the element in the region or through increasing the number of control nodes, until one of the above conditions is reached or the difference between the local minimum and maximum bound has reached some acceptable tolerance. In case of the latter, it is common to treat the element as invalid as the minimum determinant in this case is very close to zero.

\subsection{Application to bounds-preserving limiters}\label{ssec:limiter}
Another application for where finding a bound for the extrema of polynomials is of interest is in developing bounds-preserving limiters for high-order finite element methods, particularly for hyperbolic conservation laws. In these systems, it is often necessary for the solution to reside within some set of bounds $u_{\min} \leq u_h(\mathbf{x}) \leq u_{\max}$ (e.g., local maximum principle for linear transport, positivity of density in gas dynamics, etc.), which is typically enforced by applying some sort of limiting to the solution. However, applying limiting at discrete nodal points is problematic in applications requiring solution evaluation at new points, such as adaptive mesh refinement, multi-physics coupling with independent meshes/solvers, arbitrary Lagrangian--Eulerian methods, and overset meshes. Here, it is necessary for the solution to abide by these bounds across the entire solution polynomial, which is typically accomplished through limiting on Bernstein representations~\citep{Anderson2017, Glaubitz2019,Hajduk2021} or nonlinear optimization-based approaches~\citep{Dzanic2024,Dzanic2024b}.

Given that the proposed approach can yield guaranteed bounds on the high-order solution, these bounds can also be directly used to construct continuously bounds-preserving limiting approaches. In particular, we consider a discontinuous Galerkin (DG)-type approximation~\citep{Reed1973} of hyperbolic conservation laws of the form  
\begin{equation}
    \partial_t u + c{\cdot}\nabla u = 0,
\end{equation}
where $c$ is some constant advection velocity. For DG approximations of these conservation laws, it has been shown that the element-wise mean, defined as for an arbitrary element $\Omega$ as
\begin{equation}
    \overline{u} = \frac{\int_{\Omega} u \ \mathrm{d}\mathbf{x}}{\int_{\Omega} \mathrm{d}\mathbf{x}},
\end{equation}
preserves maximum principle bounds of the form
\begin{equation}
    a\leq \overline{u} (\mathbf{x}, t + \Delta t) \leq b, \quad a = \underset{\mathbf{x}}{\min}\ u(\mathbf{x}, t), \quad b =\underset{\mathbf{x}}{\max}\ u(\mathbf{x}, t),
\end{equation}
under some relatively minor assumptions on the numerical scheme and time step $\Delta t$ (for further details, we refer the reader to a series of works originating from \citet{Zhang2010}). Therefore, we can compute a limited solution $\widetilde{u_h}(\mathbf{x})$ using the squeeze-type limiter of \citet{Zhang2010} by blending the high-order solution $u_h(\mathbf{x})$ within an element $\Omega$ with its element-wise mean $\overline{u}$ as
\begin{equation}
    \widetilde{u_h}(\mathbf{x}) = \alpha u_h(\mathbf{x}) + (1 - \alpha)\overline{u},
\end{equation}
where $\alpha \in [0,1]$ is some element-wise constant convex blending coefficient. 

It can easily be shown that if $\alpha$ is computed as
\begin{equation}
    \alpha = \min \left[1,  \frac{a - \overline{u}}{u_{\min} - \overline{u}}, \frac{b - \overline{u}}{u_{\max} - \overline{u}} \right]
\end{equation}
and $u_{\min}$ and $u_{\max}$ bound the minimum and maximum of $u_h(\mathbf{x})$ within the element (i.e., $u_{\min} \leq u_h(\mathbf{x}) \leq u_{\max}\ \forall\ \mathbf{x} \in \Omega$), then the limited solution will preserve bounds across the entire element (i.e., $a \leq \widetilde{u_h}(\mathbf{x}) \leq b\ \forall\ \mathbf{x} \in \Omega$). As such, we propose to use the presented bounding technique to compute $u_{\min}$ and $u_{\max}$. Note here that while we focus on squeeze-type limiters for the linear transport equation for the purposes of this work, the methods to be presented can be applied to different limiting techniques and conservation laws which simply require bounds for the extrema of the solution. 

\subsection{Optimization process}\label{ssec:opt}
The optimization process relies on solving a constrained optimization problem for computing the optimal bounding functions. We start first with the simpler example of computing the bounding functions for a fixed set of control nodes $\boldsymbol{\eta}$. For brevity, we present this with respect to an upper bounding function, represented by its control values $\mathbf{q}$, for an arbitrary basis function $\phi(\mathbf{x})$. The problem of just checking the constraint $g(\mathbf{q}, \mathbf{x})$ in \cref{eq:obj1} is an optimization problem in itself, requiring calculating the minimum of a high-order polynomial which can be computationally expensive to compute many times. Therefore, consider instead a discrete version of the constraint, where we check the constraint at a set of $n$ equispaced sampling points as
\begin{equation}
    g(\mathbf{q}, \mathbf{x}_i) \geq 0  \ \forall \ i \in \lbrace 1, \ldots, n \rbrace.
\end{equation}

If we choose $n$ to be large, the discrete formulation quickly converges to the continuous formulation. Therefore, we can apply the optimization process to the much simpler discrete objective function 
\begin{equation}
    f(\mathbf{q}) \approx \sqrt{\frac{1}{n} \sum_{i=1}^{n} \left( L_{\boldsymbol{\eta}, \mathbf{q}}(\mathbf{x}_i) - \phi(\mathbf{x}_i) \right)^2 + \left( U_{\boldsymbol{\eta}, \mathbf{q}}(\mathbf{x}_i) - \phi(\mathbf{x}_i) \right)^2}
\end{equation}
with the above discrete constraints to obtain a \emph{candidate} solution $\mathbf{q}'$. We assume here that $n$ is chosen sufficiently large such that the \emph{exact} solution $\mathbf{q}^*$ does not differ much from the candidate solution. However, this does not ensure that the candidate solution actually satisfies the continuous constraint in \cref{eq:obj1} since the basis function may exceed the bounds outside of the sampling nodes. Therefore, we offset the control node values to account for this as 
\begin{equation}
    \mathbf{q}^* = \mathbf{q}' + \Delta \mathbf{q} + \epsilon,
\end{equation}
where 
\begin{equation}
    \Delta \mathbf{q} =  -\underset{\mathbf{x}}\min \left (0, L_{\boldsymbol{\eta}, \mathbf{q}'}(\mathbf{x}) - \phi(\mathbf{x}) \right) = - \underset{\mathbf{x}}\min\ (0, g(\mathbf{q}', \mathbf{x}))
\end{equation}
is the maximum amount the continuous formulation violates the constraint functional and $\epsilon = 10^{-6}$ is a small numerical tolerance. With this approach, the discrete formulation can be optimized over many sampling points relatively quickly while the problem of calculating the minimum of a high-order polynomial is required only once at the end. This greatly reduces the computational complexity of the problem, and the differences in the end results are essentially negligible for large $n$ (in our case, $n = 1000$). We perform this constrained optimization using the sequential least-squares quadratic programming (SLSQP) algorithm in the SciPy package. Furthermore, symmetry in the bounding boxes with respect to symmetric basis functions was enforced by optimizing on one ``side'' and reflecting the boxes, and, for higher-dimensional bases, we compute the bounding boxes as tensor-products of the one-dimensional bounding boxes.

The optimization problem for finding the optimal control nodes is significantly more complex, with each step of this optimization process requiring computing the bounding boxes for a fixed set of nodes as above. In addition to the increased computational cost, the distribution of the nodes themselves have constraints which must be enforced, namely that they must be symmetric about the origin, must have nodes on the endpoints $x = \pm 1$, and must be bounded by the element domain $\Omega = [-1,1]$ and distributed in increasing order. The first two are simple to enforce by ensuring only some nodes are free to move -- for example, only nodes on the left side (excluding the endpoint and midpoint if applicable) are part of the optimization process while the right side is taken as a reflection of the left. For the last part, this can be imposed by a further set of constraints on the optimizer (e.g., $q_{i+1} - q_{i} \geq 0$, $  q_i + 1  \geq 0$, $1 - q_i  \geq 0$), but we found that with this many constraints, optimizers often struggled to converge to the optimal solution. 

Instead, we convert a (subset of) this constrained optimization problem into an unconstrained optimization problem via an auxiliary variable-type approach. We optimize instead for a set of $M-1$ auxiliary variables $\mathbf{z}$, which are unbounded, and transform them to the control nodes as
\begin{equation}
    q_i = -1 + 2\frac{\sum_{j=1}^{i} \exp (z_j)}{\sum_{j=1}^{M-1} \exp (z_j)}.
\end{equation}
As it can be seen, this guarantees that $q_{i+1} > q_{i}$, $\min {\mathbf{q}} = -1$, $\max {\mathbf{q}} = 1$, and $ -1 \leq q_i \leq 1$ without imposing any further constraints on the optimizer at the expense of one additional variable. Note that the symmetry arguments from before can also be applied here to reduce the number of variables that need to be optimized over. We similarly perform this optimization using the Broyden--Fletcher--Goldfarb--Shanno (BFGS) algorithm with basin-hopping in the SciPy package, with the above SLSQP optimizer for the inner optimization steps for computing the bounding boxes for the current set of control nodes. The implementation of this optimization process is included as Python code in the electronic supplementary material, and some tabulated examples of the bounding box control nodes/values are presented in \cref{app:bounds}.

\subsection{Overview}
We present here a brief overview of the approach as applied to computing bounds on a high-order polynomial solution within an element. We assume here the optimal bounding boxes for the chosen basis functions have been precomputed as per \cref{ssec:opt} and stored for loading. Then, for each element:
\begin{enumerate}
    \item Load two tables of $N \times M$ entries for the given basis, consisting of $N$ basis functions and $M$ optimal control node values. Denote $q^{ij}_-$/$q^{ij}_+$ as the lower/upper control node values for basis function $\phi_i$ and control node $\eta_j$.
    \item Compute projection coefficients (from \cref{eq:projcoeff}) via quadrature. Subtract linear basis $u_{LO}(\mathbf{x})$ from solution $u_h(\mathbf{x})$ to compute high-order fluctuations $u_h'(\mathbf{x})$.
    \item Loop over control nodes $1 \leq i \leq M$:
      \begin{enumerate}
        \item Compute fluctuation bounds at control node: 
        \begin{equation*}
            u'_{\min,i} = \sum_{i =1}^N \min(u_i' q^{ij}_-, u_i' q^{ij}_+) \quad \mathrm{and} \quad u'_{\max,i} = \sum_{i =1}^N \max(u_i' q^{ij}_-, u_i' q^{ij}_+).
        \end{equation*}
        \item Compute solution bounds at control node: $u_{\min,i} = u'_{\min,i} + u_{LO}(\mathbf{x}_i)$ and $u_{\max,i} = u'_{\max,i} + u_{LO}(\mathbf{x}_i)$
      \end{enumerate}
    \item If desired (for increased accuracy), increase $M$ and repeat or interpolate solution onto sub-element and repeat. 
\end{enumerate}

\section{Results}\label{sec:results}
\begin{figure}[htbp!]
        \centering
        \subfloat[$N=3$]{\adjustbox{width=0.32\linewidth, valign=b}{\begin{tikzpicture}[spy using outlines={rectangle, height=3cm,width=2.5cm, magnification=3, connect spies},
font=\large]
    \begin{semilogyaxis}
    [
        axis line style={latex-latex},
        axis y line=left,
        axis x line=left,
        clip mode=individual,
        xlabel = {$M$},
        ylabel = {$\varepsilon_2$},
        xmin = 3, xmax = 20,
        ymin = 1e-3, ymax = 2e-1,
        legend cell align={left},
        legend style={font=\normalsize, at={(0.97, 0.97)}, anchor=north east},
        x tick label style={/pgf/number format/.cd, fixed, fixed zerofill, precision=0, /tikz/.cd},
        y tick label style={/pgf/number format/.cd, fixed, precision=2, /tikz/.cd},	
        grid=both,
        grid style={line width=.1pt, draw=gray!10},
        major grid style={line width=.2pt,draw=gray!50},
        minor x tick num=4,
        minor y tick num=4,
    ]
        
        \addplot[color=PlotColor1, style={thick}, mark=*, mark options={scale=0.5}] table[x=M, y expr={\thisrow{gl}/3}, col sep=comma]{./figs/data/basis_l2errs_N3.csv};
        \addlegendentry{Gauss--Legendre};
        
        \addplot[color=PlotColor2, style={thick}, mark=*, mark options={scale=0.5}] table[x=M, y expr={\thisrow{gll}/3}, col sep=comma]{./figs/data/basis_l2errs_N3.csv};
        \addlegendentry{Gauss--Lobatto};
        
        \addplot[color=PlotColor3, style={thick}, mark=*, mark options={scale=0.5}] table[x=M, y expr={\thisrow{cheb}/3}, col sep=comma]{./figs/data/basis_l2errs_N3.csv};
        \addlegendentry{Chebyshev};
        
        \addplot[color=PlotColor4, style={thick}, mark=*, mark options={scale=0.5}] table[x=M, y expr={\thisrow{eq}/3}, col sep=comma]{./figs/data/basis_l2errs_N3.csv};
        \addlegendentry{Equispaced};
        
        \addplot[color=PlotColor5, style={thick}, mark=*, mark options={scale=0.5}] table[x=M, y expr={\thisrow{opt}/3}, col sep=comma]{./figs/data/basis_l2errs_N3.csv};
        \addlegendentry{Optimized};
    
    \end{semilogyaxis}
\end{tikzpicture}}}
        \subfloat[$N=4$]{\adjustbox{width=0.32\linewidth, valign=b}{\begin{tikzpicture}[spy using outlines={rectangle, height=3cm,width=2.5cm, magnification=3, connect spies},
font=\large]
    \begin{semilogyaxis}
    [
        axis line style={latex-latex},
        axis y line=left,
        axis x line=left,
        clip mode=individual,
        xlabel = {$M$},
        ylabel = {$\varepsilon_2$},
        xmin = 3, xmax = 20,
        ymin = 3e-3, ymax =3e-1,
        legend cell align={left},
        legend style={font=\scriptsize, at={(0.97, 0.97)}, anchor=north east},
        x tick label style={/pgf/number format/.cd, fixed, fixed zerofill, precision=0, /tikz/.cd},
        y tick label style={/pgf/number format/.cd, fixed, precision=2, /tikz/.cd},	
        grid=both,
        grid style={line width=.1pt, draw=gray!10},
        major grid style={line width=.2pt,draw=gray!50},
        minor x tick num=4,
        minor y tick num=4,
    ]
        
        \addplot[color=PlotColor1, style={thick}, mark=*, mark options={scale=0.5}] table[x=M, y expr={\thisrow{gl}/4}, col sep=comma]{./figs/data/basis_l2errs_N4.csv};
        
        \addplot[color=PlotColor2, style={thick}, mark=*, mark options={scale=0.5}] table[x=M, y expr={\thisrow{gll}/4}, col sep=comma]{./figs/data/basis_l2errs_N4.csv};
        
        \addplot[color=PlotColor3, style={thick}, mark=*, mark options={scale=0.5}] table[x=M, y expr={\thisrow{cheb}/4}, col sep=comma]{./figs/data/basis_l2errs_N4.csv};
        
        \addplot[color=PlotColor4, style={thick}, mark=*, mark options={scale=0.5}] table[x=M, y expr={\thisrow{eq}/4}, col sep=comma]{./figs/data/basis_l2errs_N4.csv};
        
        \addplot[color=PlotColor5, style={thick}, mark=*, mark options={scale=0.5}] table[x=M, y expr={\thisrow{opt}/4}, col sep=comma]{./figs/data/basis_l2errs_N4.csv};
    
    \end{semilogyaxis}
\end{tikzpicture}}}
        \subfloat[$N=5$]{\adjustbox{width=0.32\linewidth, valign=b}{\begin{tikzpicture}[spy using outlines={rectangle, height=3cm,width=2.5cm, magnification=3, connect spies},
font=\large]
    \begin{semilogyaxis}
    [
        axis line style={latex-latex},
        axis y line=left,
        axis x line=left,
        clip mode=individual,
        xlabel = {$M$},
        ylabel = {$\varepsilon_2$},
        xmin = 3, xmax = 20,
        ymin = 5e-3, ymax = 4e-1,
        legend cell align={left},
        legend style={font=\scriptsize, at={(0.97, 0.97)}, anchor=north east},
        x tick label style={/pgf/number format/.cd, fixed, fixed zerofill, precision=0, /tikz/.cd},
        y tick label style={/pgf/number format/.cd, fixed, precision=2, /tikz/.cd},	
        grid=both,
        grid style={line width=.1pt, draw=gray!10},
        major grid style={line width=.2pt,draw=gray!50},
        minor x tick num=4,
        minor y tick num=4,
    ]
        
        \addplot[color=PlotColor1, style={thick}, mark=*, mark options={scale=0.5}] table[x=M, y expr={\thisrow{gl}/5}, col sep=comma]{./figs/data/basis_l2errs_N5.csv};
        
        \addplot[color=PlotColor2, style={thick}, mark=*, mark options={scale=0.5}] table[x=M, y expr={\thisrow{gll}/5}, col sep=comma]{./figs/data/basis_l2errs_N5.csv};
        
        \addplot[color=PlotColor3, style={thick}, mark=*, mark options={scale=0.5}] table[x=M, y expr={\thisrow{cheb}/5}, col sep=comma]{./figs/data/basis_l2errs_N5.csv};
        
        \addplot[color=PlotColor4, style={thick}, mark=*, mark options={scale=0.5}] table[x=M, y expr={\thisrow{eq}/5}, col sep=comma]{./figs/data/basis_l2errs_N5.csv};
        
        \addplot[color=PlotColor5, style={thick}, mark=*, mark options={scale=0.5}] table[x=M, y expr={\thisrow{opt}/5}, col sep=comma]{./figs/data/basis_l2errs_N5.csv};
    
    \end{semilogyaxis}
\end{tikzpicture}}}
        \newline
        \subfloat[$N=6$]{\adjustbox{width=0.32\linewidth, valign=b}{\begin{tikzpicture}[spy using outlines={rectangle, height=3cm,width=2.5cm, magnification=3, connect spies},
font=\large]
    \begin{semilogyaxis}
    [
        axis line style={latex-latex},
        axis y line=left,
        axis x line=left,
        clip mode=individual,
        xlabel = {$M$},
        ylabel = {$\varepsilon_2$},
        xmin = 3, xmax = 20,
        ymin = 1e-2, ymax = 5e-1,
        legend cell align={left},
        legend style={font=\scriptsize, at={(0.97, 0.97)}, anchor=north east},
        x tick label style={/pgf/number format/.cd, fixed, fixed zerofill, precision=0, /tikz/.cd},
        y tick label style={/pgf/number format/.cd, fixed, precision=2, /tikz/.cd},	
        grid=both,
        grid style={line width=.1pt, draw=gray!10},
        major grid style={line width=.2pt,draw=gray!50},
        minor x tick num=4,
        minor y tick num=4,
    ]
        
        \addplot[color=PlotColor1, style={thick}, mark=*, mark options={scale=0.5}] table[x=M, y expr={\thisrow{gl}/6}, col sep=comma]{./figs/data/basis_l2errs_N6.csv};
        
        \addplot[color=PlotColor2, style={thick}, mark=*, mark options={scale=0.5}] table[x=M, y expr={\thisrow{gll}/6}, col sep=comma]{./figs/data/basis_l2errs_N6.csv};
        
        \addplot[color=PlotColor3, style={thick}, mark=*, mark options={scale=0.5}] table[x=M, y expr={\thisrow{cheb}/6}, col sep=comma]{./figs/data/basis_l2errs_N6.csv};
        
        \addplot[color=PlotColor4, style={thick}, mark=*, mark options={scale=0.5}] table[x=M, y expr={\thisrow{eq}/6}, col sep=comma]{./figs/data/basis_l2errs_N6.csv};
        
        \addplot[color=PlotColor5, style={thick}, mark=*, mark options={scale=0.5}] table[x=M, y expr={\thisrow{opt}/6}, col sep=comma]{./figs/data/basis_l2errs_N6.csv};
    
    \end{semilogyaxis}
\end{tikzpicture}}}
        \subfloat[$N=7$]{\adjustbox{width=0.32\linewidth, valign=b}{\begin{tikzpicture}[spy using outlines={rectangle, height=3cm,width=2.5cm, magnification=3, connect spies},
font=\large]
    \begin{semilogyaxis}
    [
        axis line style={latex-latex},
        axis y line=left,
        axis x line=left,
        clip mode=individual,
        xlabel = {$M$},
        ylabel = {$\varepsilon_2$},
        xmin = 3, xmax = 20,
        ymin = 1e-2, ymax = 6e-1,
        legend cell align={left},
        legend style={font=\scriptsize, at={(0.97, 0.97)}, anchor=north east},
        x tick label style={/pgf/number format/.cd, fixed, fixed zerofill, precision=0, /tikz/.cd},
        y tick label style={/pgf/number format/.cd, fixed, precision=2, /tikz/.cd},	
        grid=both,
        grid style={line width=.1pt, draw=gray!10},
        major grid style={line width=.2pt,draw=gray!50},
        minor x tick num=4,
        minor y tick num=4,
    ]
        
        \addplot[color=PlotColor1, style={thick}, mark=*, mark options={scale=0.5}] table[x=M, y expr={\thisrow{gl}/7}, col sep=comma]{./figs/data/basis_l2errs_N7.csv};
        
        \addplot[color=PlotColor2, style={thick}, mark=*, mark options={scale=0.5}] table[x=M, y expr={\thisrow{gll}/7}, col sep=comma]{./figs/data/basis_l2errs_N7.csv};
        
        \addplot[color=PlotColor3, style={thick}, mark=*, mark options={scale=0.5}] table[x=M, y expr={\thisrow{cheb}/7}, col sep=comma]{./figs/data/basis_l2errs_N7.csv};
        
        \addplot[color=PlotColor4, style={thick}, mark=*, mark options={scale=0.5}] table[x=M, y expr={\thisrow{eq}/7}, col sep=comma]{./figs/data/basis_l2errs_N7.csv};
        
        \addplot[color=PlotColor5, style={thick}, mark=*, mark options={scale=0.5}] table[x=M, y expr={\thisrow{opt}/7}, col sep=comma]{./figs/data/basis_l2errs_N7.csv};
    
    \end{semilogyaxis}
\end{tikzpicture}}}
        \subfloat[$N=8$]{\adjustbox{width=0.32\linewidth, valign=b}{\begin{tikzpicture}[spy using outlines={rectangle, height=3cm,width=2.5cm, magnification=3, connect spies},
font=\large]
    \begin{semilogyaxis}
    [
        axis line style={latex-latex},
        axis y line=left,
        axis x line=left,
        clip mode=individual,
        xlabel = {$M$},
        ylabel = {$\varepsilon_2$},
        xmin = 3, xmax = 20,
        ymin = 1e-2, ymax = 6e-1,
        legend cell align={left},
        legend style={font=\scriptsize, at={(0.97, 0.97)}, anchor=north east},
        x tick label style={/pgf/number format/.cd, fixed, fixed zerofill, precision=0, /tikz/.cd},
        y tick label style={/pgf/number format/.cd, fixed, precision=2, /tikz/.cd},	
        grid=both,
        grid style={line width=.1pt, draw=gray!10},
        major grid style={line width=.2pt,draw=gray!50},
        minor x tick num=4,
        minor y tick num=4,
    ]
        
        \addplot[color=PlotColor1, style={thick}, mark=*, mark options={scale=0.5}] table[x=M, y expr={\thisrow{gl}/8}, col sep=comma]{./figs/data/basis_l2errs_N8.csv};
        
        \addplot[color=PlotColor2, style={thick}, mark=*, mark options={scale=0.5}] table[x=M, y expr={\thisrow{gll}/8}, col sep=comma]{./figs/data/basis_l2errs_N8.csv};
        
        \addplot[color=PlotColor3, style={thick}, mark=*, mark options={scale=0.5}] table[x=M, y expr={\thisrow{cheb}/8}, col sep=comma]{./figs/data/basis_l2errs_N8.csv};
        
        \addplot[color=PlotColor4, style={thick}, mark=*, mark options={scale=0.5}] table[x=M, y expr={\thisrow{eq}/8}, col sep=comma]{./figs/data/basis_l2errs_N8.csv};
        
        \addplot[color=PlotColor5, style={thick}, mark=*, mark options={scale=0.5}] table[x=M, y expr={\thisrow{opt}/8}, col sep=comma]{./figs/data/basis_l2errs_N8.csv};
    
    \end{semilogyaxis}
\end{tikzpicture}}}
        \newline
        \caption{Average $L^2$ error of the bounding boxes for the $N-1$ order Gauss--Lobatto nodal interpolating basis functions as a function of the control node placements and number of control nodes $M$. }
        \label{fig:basis_l2errs}  
    \end{figure}
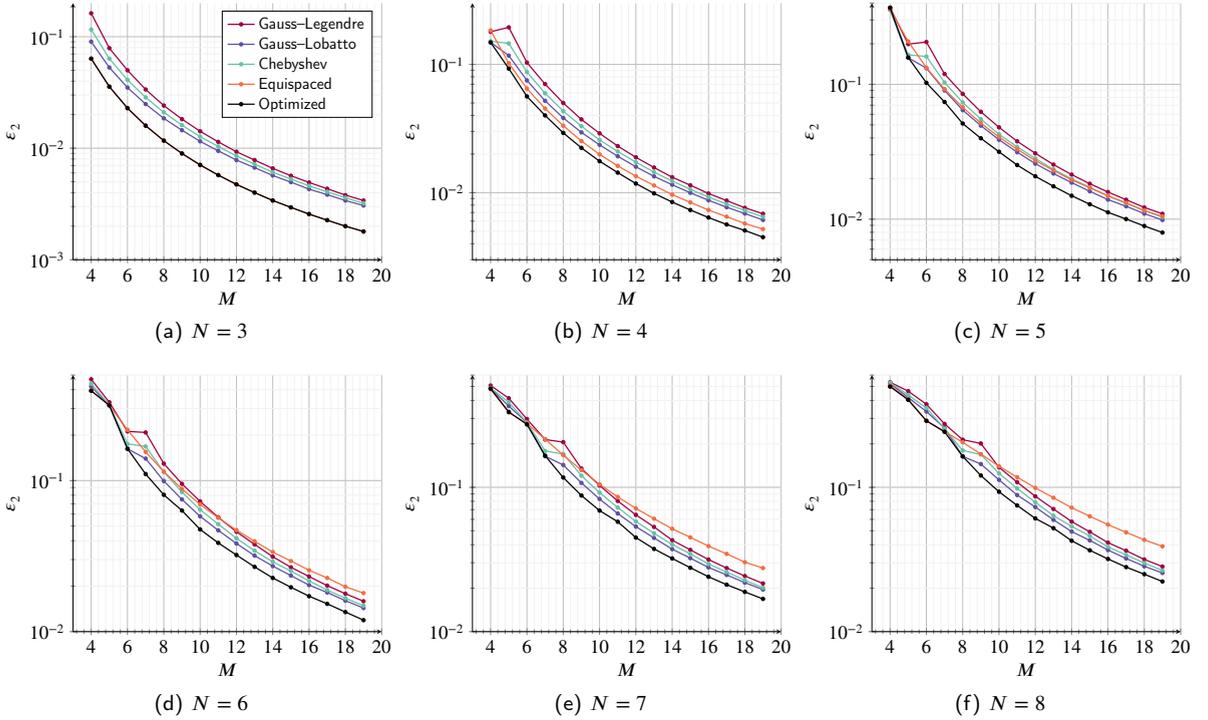

We first look at the efficacy of the proposed approach in terms of bounding the basis functions themselves for a finite element approximation with $N$ solution nodes and $M$ control nodes. The finite element basis functions were taken as the $N-1$ nodal interpolating functions on the $N$ Gauss--Lobatto nodes for a one-dimensional element. For each $N$, the error was computed with respect to the optimized bounding boxes computed with $M$ control nodes placed on the Gauss--Legendre nodes (consisting of the $M-2$ Gauss-Legendre nodes and the endpoints), the Gauss--Lobatto nodes, the Chebyshev nodes, equispaced nodes, and optimized nodes as proposed in \cref{ssec:bbox}. The error, denoted by $\varepsilon_2$, was computed with respect to the sum of the $L^2$ norm of the bounding box errors (i.e., $\varepsilon_2 = f(\boldsymbol{\eta}, \mathbf{q})$ in \cref{eq:obj2}). The errors at varying values of $N$ and $M$ are shown in \cref{fig:basis_l2errs} for the different control node points. It can be seen that, as expected, the optimized nodes result in the lowest average error. At low $N$, the equispaced nodes were typically the second-best control node set, with $N=3$ showing essentially identical error between the equispaced nodes and optimized nodes. However, as $N$ increased, the equispaced node set quickly became suboptimal, and the performance of the Gauss--Lobatto node set was closer to the optimal nodes, followed by the Chebyshev and Gauss--Legendre nodes. For all node sets, the error behaved approximately as $\mathcal O(M^{-2})$, indicating second-order convergence with respect to average control node spacing which is consistent with the property that a piecewise linear interpolant of a smooth function converges with second-order accuracy in the $L^2$ norm. We note here that the error between the best and worst performing node sets typically only differed by a factor of 2-3.

\subsection{High-order mesh validity checks}
\begin{figure}[htbp!]
        \centering
        \adjustbox{width=0.48\linewidth, valign=b}{\begin{tikzpicture}[spy using outlines={rectangle, height=3cm,width=2.5cm, magnification=3, connect spies}]
    \begin{semilogyaxis}
    [
        axis line style={latex-latex},
        axis y line=left,
        axis x line=left,
        clip mode=individual,
        xlabel = {Subdivision levels},
        ylabel = {$\text{min} (|\mathbf{J}|)$ estimate},
        xmin = 0, xmax = 8,
        ymin = 1e-4, ymax = 1,
        ytick = {1e-4, 1e-3, 1e-2, 1e-1, 1},
        yticklabels = {$-10^{-4}$, $-10^{-3}$, $-10^{-2}$, $-10^{-1}$, $-1$},
        y dir=reverse,
        legend cell align={left},
        legend style={font=\scriptsize, at={(0.97, 0.03)}, anchor=south east},
        x tick label style={/pgf/number format/.cd, fixed, fixed zerofill, precision=0, /tikz/.cd},
        y tick label style={/pgf/number format/.cd, fixed, precision=2, /tikz/.cd},
        grid=both,
        grid style={line width=.1pt, draw=gray!10},
        major grid style={line width=.2pt,draw=gray!50},
        minor x tick num=1,
    ]
        
        \addplot[color=PlotColor1, style={thick}, mark=*, mark options={scale=0.5}] coordinates {
        (0, 4.1912e-01)
        (1, 1.0871e-01)
        (2, 9.1725e-03)
        (3, 4.2475e-03) 
        (4, 1.7850e-03)
        (5, 5.5382e-04)
        (6, 2.3632e-04)
        (7, 2.2568e-04)
    };
        \addlegendentry{Bernstein};

        \addplot[color=PlotColor2, style={thick}, mark=*, mark options={scale=0.5}] coordinates {
        (0,1.99E-01)
        (1,2.10E-02)
        (2,1.46E-03)
        (3,4.65E-04)
        (4,2.82E-04)
        
    };
        \addlegendentry{Present work ($M = N$)};

        \addplot[color=PlotColor3, style={thick}, mark=*, mark options={scale=0.5}] coordinates {
        (0,9.73E-02)
        (1,8.83E-03)
        (2,7.89E-04)
        (3,2.58E-04)
    };
        \addlegendentry{Present work ($M = N+1$)};

        \addplot[color=PlotColor4, style={thick}, mark=*, mark options={scale=0.5}] coordinates {
        (0,3.10E-02)
        (1,3.87E-03)
        (2,3.22E-04)
    };
        \addlegendentry{Present work ($M = N+2$)};
        
        \addplot[color=black, style={very thick, densely dotted}] coordinates {
        (0,0.000215656)
        (8,0.000215656)
    };
        \addlegendentry{Exact};




    \end{semilogyaxis}
\end{tikzpicture}}
        \caption{Estimate of the minimum determinant of the element transformation Jacobian with respect to subdivision levels for a second-order quadrilateral element as computed with the proposed approach and the Bernstein approach.}
        \label{fig:quad_detj_comparison}  
    \end{figure}
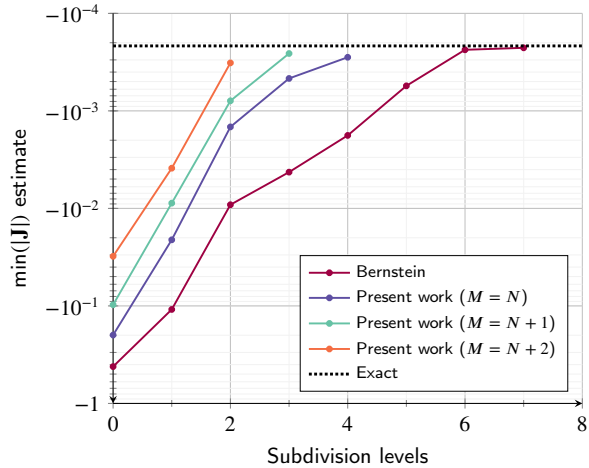

The proposed approach was then applied to check the validity (i.e., the positivity of the determinant of the element transformation Jacobian) of high-order curved meshes. We focus here on tensor-product elements, specifically quadrilateral elements, but the general techniques presented in this work can extend to any element type. For all of the numerical experiments, we use the modular open-source C++ FEM library MFEM~\citep{mfem,MFEM2024}, and some of the methods presented are available as examples in the MFEM package. As a first example, we consider the test of simply bounding $\min(|J|)$ within a second-order quadrilateral element (i.e., $\mathbb P_2$ in maximal order), for which the Jacobian determinant is $\mathbb P_3$ in maximal order (i.e., $N=4$ in each dimension). The transformation was chosen such that the element was inverted in a small region, with $\min(|J|) = -0.0002156$. The minimum bound with respect to number of subdivision levels is shown in \cref{fig:quad_detj_comparison} as computed by the proposed approach (with $M = N$, $N+1$, and $N+2$) and the Bernstein approach. It can be seen that the proposed approach yields much tighter bounds for the minimum determinant than the Bernstein approach, with around an order of magnitude improvement. As expected, $M = N+2$ yielded the tightest bounds, requiring only 2 subdivisions to reach a tolerance of $10^{-4}$ between the minimum/maximum bound, whereas the Bernstein approach required 6. With $M = N$ and $M = N+1$, 4 and 3 subdivision levels were required, respectively, which was still a significant improvement over the Bernstein approach. We remark here that one of the benefits of the proposed approach is that the bounds are local, such that the initial subdivision can be performed only on the subcell (i.e., the region between control nodes) where the bounds exceed the tolerances, instead of across the entire element which is typical of standard Bernstein-type approaches. 

\begin{figure}[htbp!]
        \centering
        \hspace{5em}
        \subfloat[Exact (brute force)]{
        \adjustbox{width=0.48\linewidth, valign=b}{\includegraphics[]{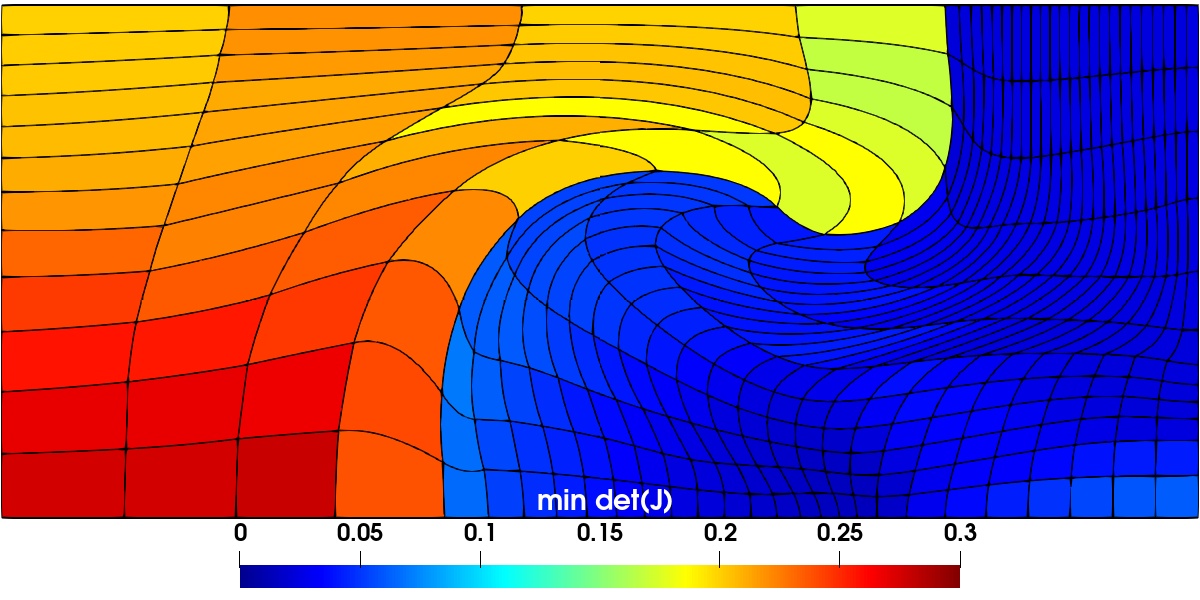}}}
        \newline
        \subfloat[Proposed approach]{
        \adjustbox{width=0.48\linewidth, valign=b}{\includegraphics[]{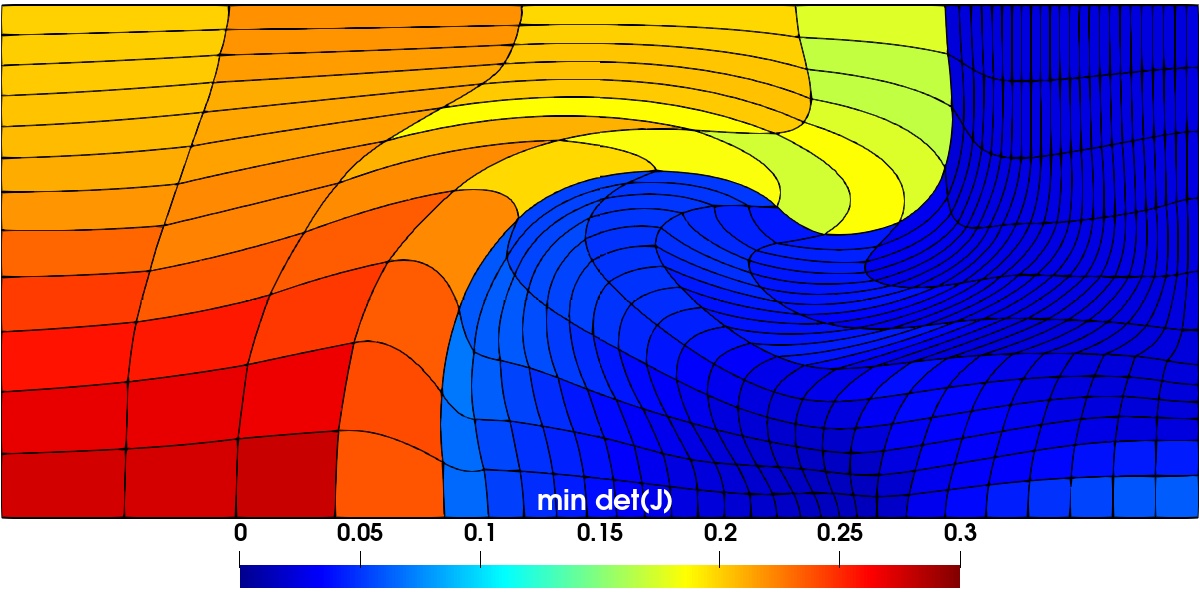}}}
        \subfloat[Bernstein approach]{
        \adjustbox{width=0.48\linewidth, valign=b}{\includegraphics[]{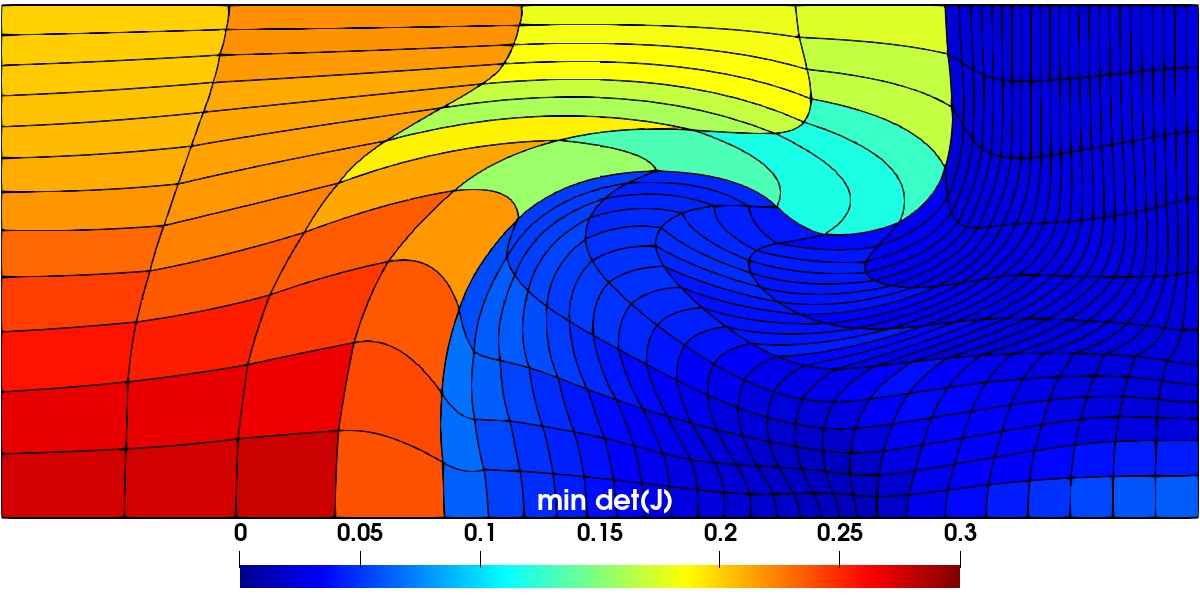}}}
        \newline
        \subfloat[Error with proposed approach]{
        \adjustbox{width=0.48\linewidth, valign=b}{\includegraphics[]{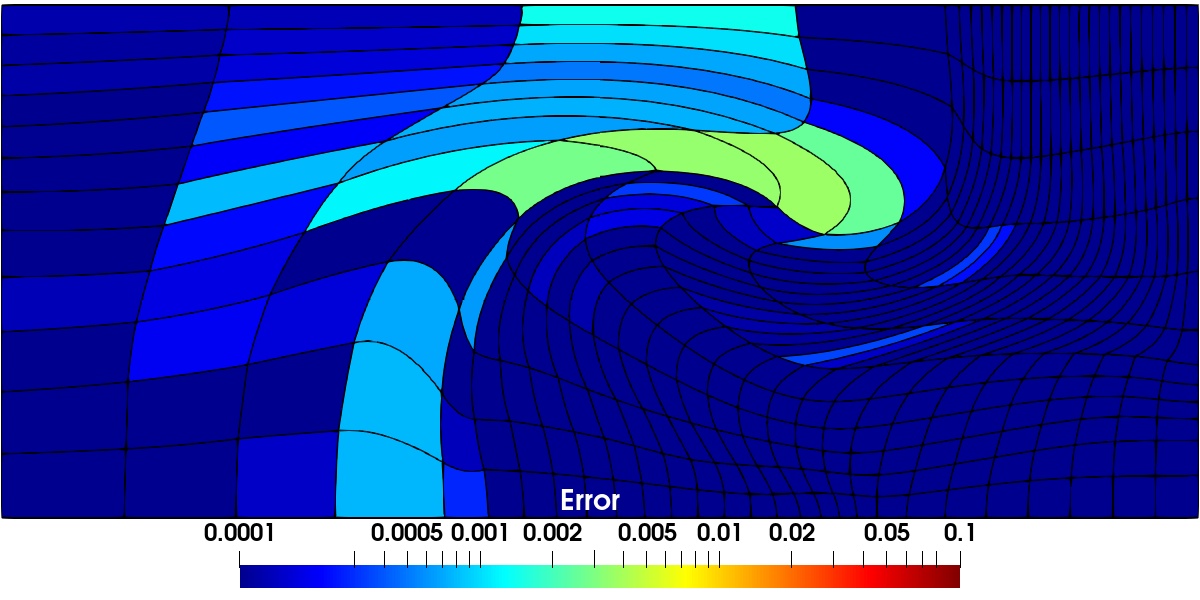}}}
        \subfloat[Error with Bernstein approach]{
        \adjustbox{width=0.48\linewidth, valign=b}{\includegraphics[]{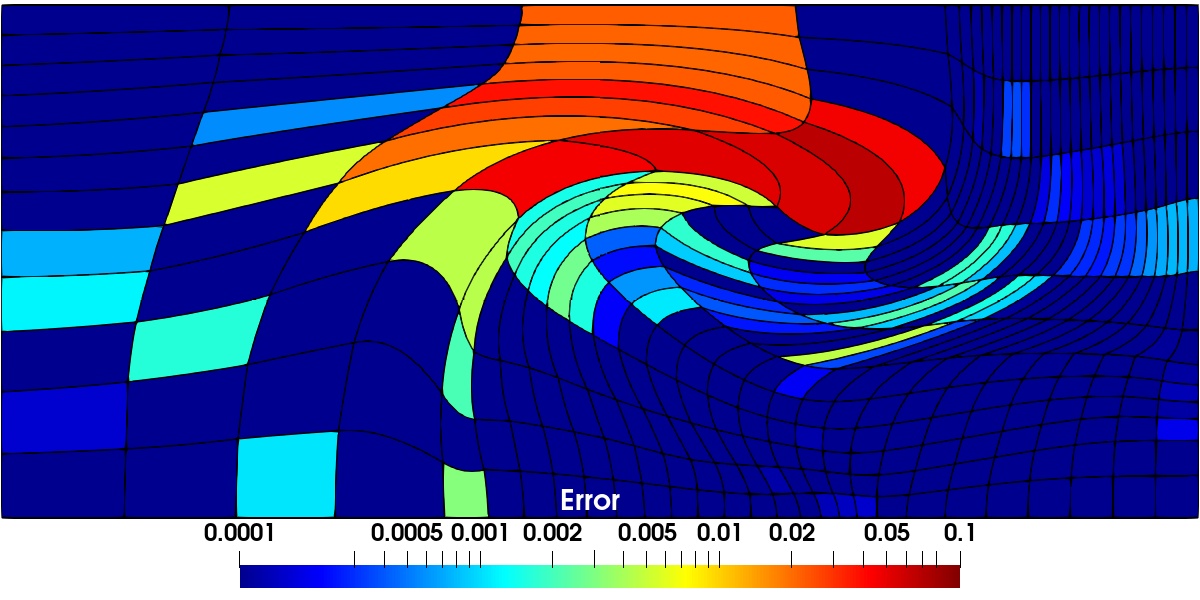}}}
        \caption{Mesh obtained from a Lagrangian hydrodynamics simulation of the triple-point shock interaction problem. Element-wise minimum mesh Jacobian determinant shown as computed by a brute force approach (top), proposed approach (middle left), and Bernstein approach (middle right). Errors for the proposed and Bernstein approaches shown on bottom row.}
        \label{fig:triplepoint_mesh_example}
    \end{figure}

We further showcase the proposed approach in the context of mesh validity checks in Lagrangian hydrodynamics. Here, the formation of non-physical elements (e.g., inverted elements) due to large deformations in high-speed flow regimes can degrade solver convergence and cause issues with solution interpolation for adaptive mesh refinement and multi-physics solvers. As such, effective techniques for ensuring the positivity of the Jacobian determinant of mesh elements can be highly beneficial in these applications, and methods which can better bound the determinant (from below) prevent the spurious flagging of valid elements which can degrade performance and accuracy. We study the mesh from a Lagrangian hydrodynamics simulation of the triple-point problem of \citet{Kucharik2010} as computed with the Laghos solver~\citep{Dobrev2012} to evaluate the proposed approach in terms of estimating the mesh Jacobian determinant. The minimum mesh Jacobian determinant in each element is shown in \cref{fig:triplepoint_mesh_example} as computed by an exact (brute force sampling) approach, the proposed bounding approach, and a Bernstein-type approach. It can be seen that the results of the proposed approach better represent the true minimum determinants in the mesh compared to the Bernstein approach, particularly so in the highly distorted region around the material interfaces where the nonlinear terms in the determinant representation are strongest. This is most evident in the error diagrams, where the proposed approach yielded maximum errors that were 1-2 orders of magnitude lower than the Bernstein approach, with relative errors on the order of only a few percent. 

As a final demonstration of the proposed approach in validating high-order meshes, we implement this approach within an $r$-adaptive mesh optimization framework to ensure mesh validity during optimization. Our mesh optimization framework is based on the Target Matrix Optimization Paradigm (TMOP)~\citep{dobrev2019target}, where node movement is driven by variational minimization of a functional that depends on the current and \emph{target} Jacobian of the transformation of each mesh element. In standard practice~\citep{dobrev2019target,aparicio2024defining}, the validity of mesh elements is checked at a set of quadrature points through a line-search procedure. Specifically, if the node displacement determined by the mesh optimization iteration results in an invalid mesh at any of the quadrature points, it is iteratively scaled down until the resulting mesh is valid. However, this results in a similar problem as with the previous example, where element validity can only be checked at discrete points which may cause solver divergence if the quadrature points used in the simulation change. With the proposed approach, this validity can be checked for the whole element, which guarantees validity for any set of quadrature points.

\begin{figure}[htbp!]
        \centering
        \subfloat[Optimization without proposed approach]{
        \adjustbox{width=0.32\linewidth, valign=b}{\includegraphics[]{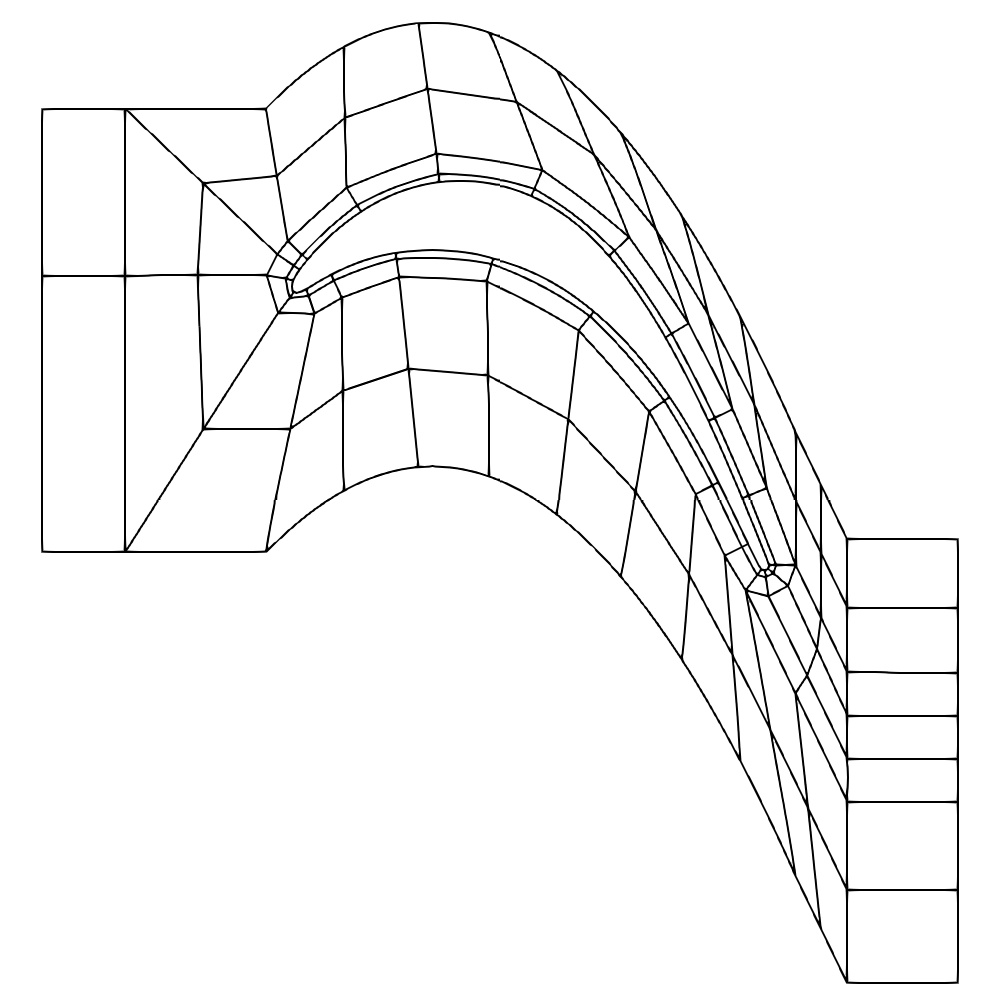}}
        \adjustbox{width=0.32\linewidth, valign=b}{\includegraphics[]{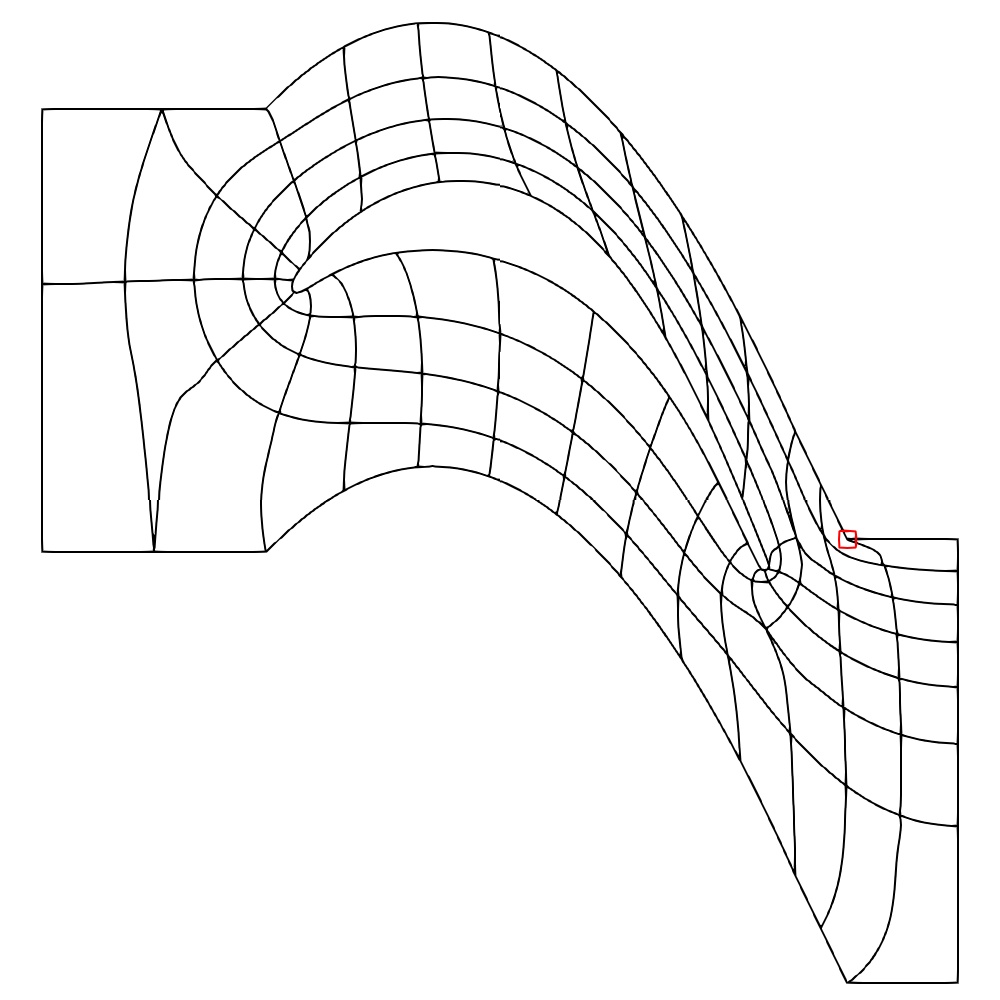}}
        \adjustbox{width=0.32\linewidth, valign=b}{\includegraphics[]{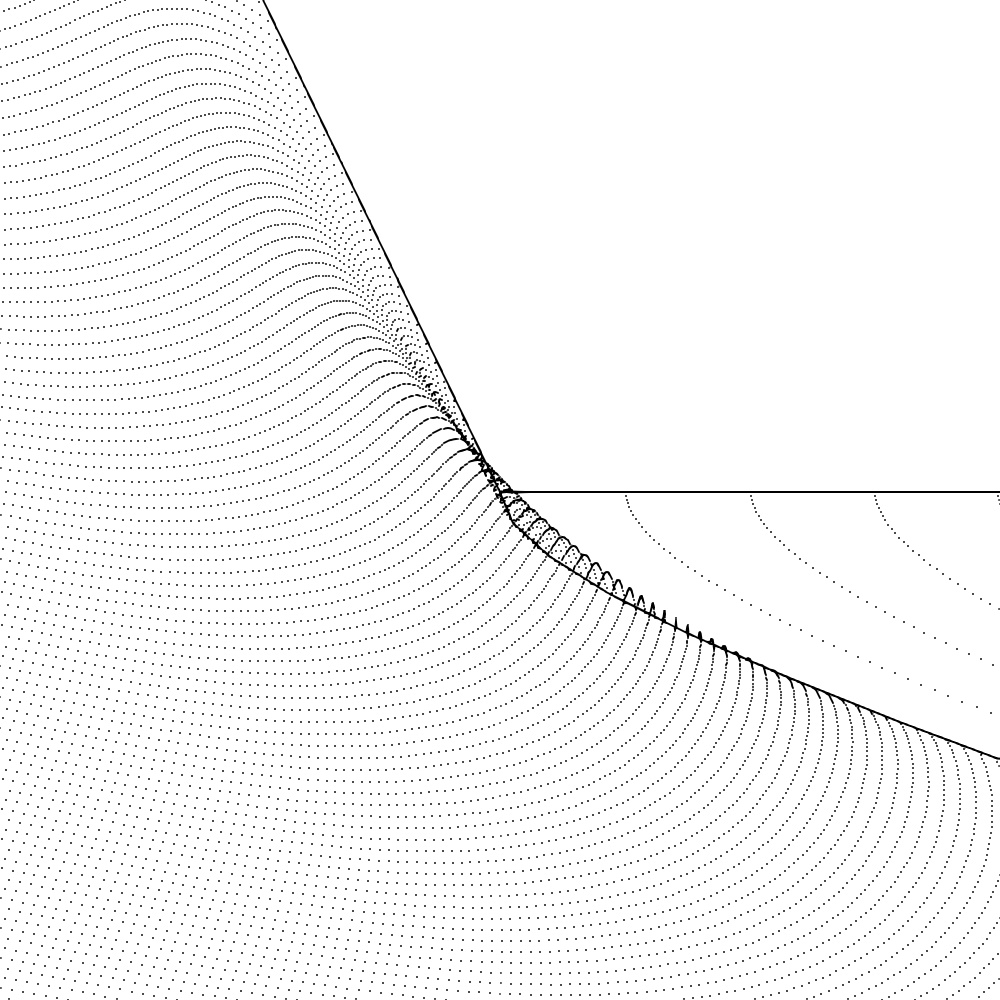}}}
        \newline
        \subfloat[Optimization with proposed approach]{
        \adjustbox{width=0.32\linewidth, valign=b}{\includegraphics[]{figs/blade_original_mesh.jpg}}
        \adjustbox{width=0.32\linewidth, valign=b}{\includegraphics[]{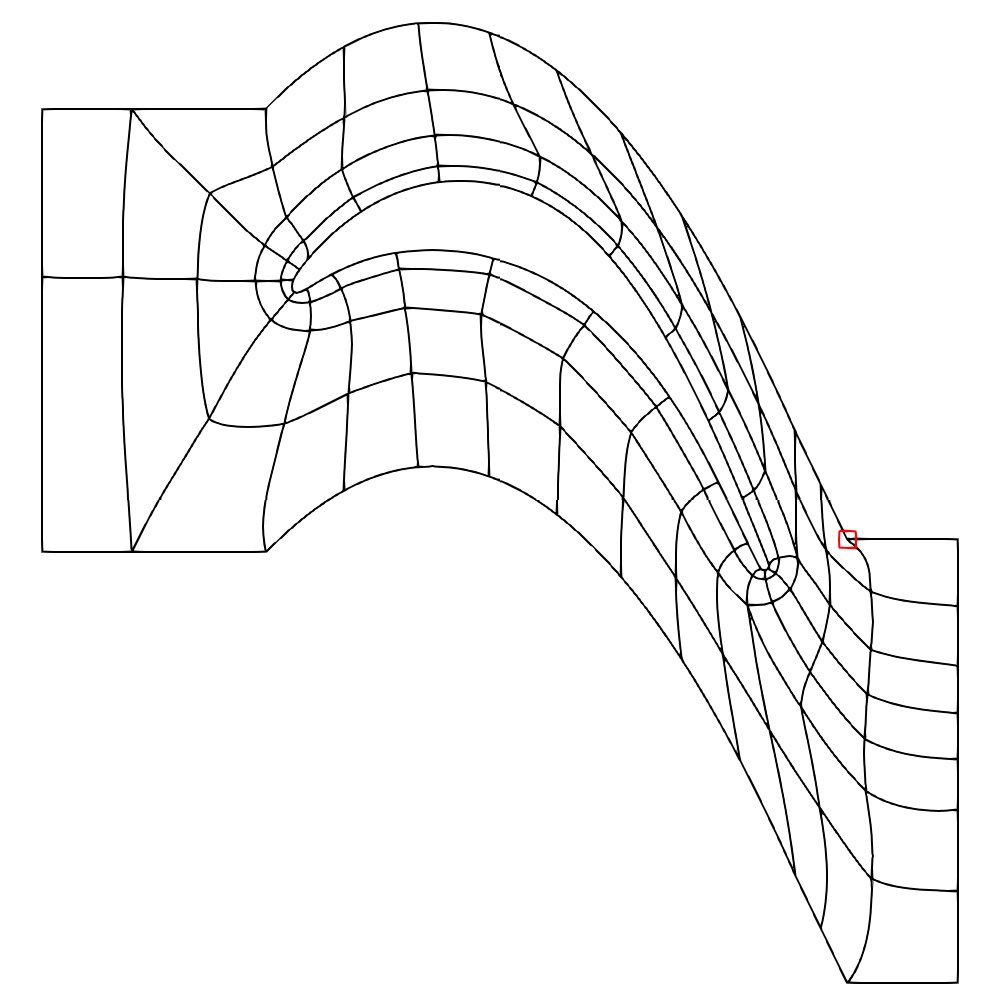}}
        \adjustbox{width=0.32\linewidth, valign=b}{\includegraphics[]{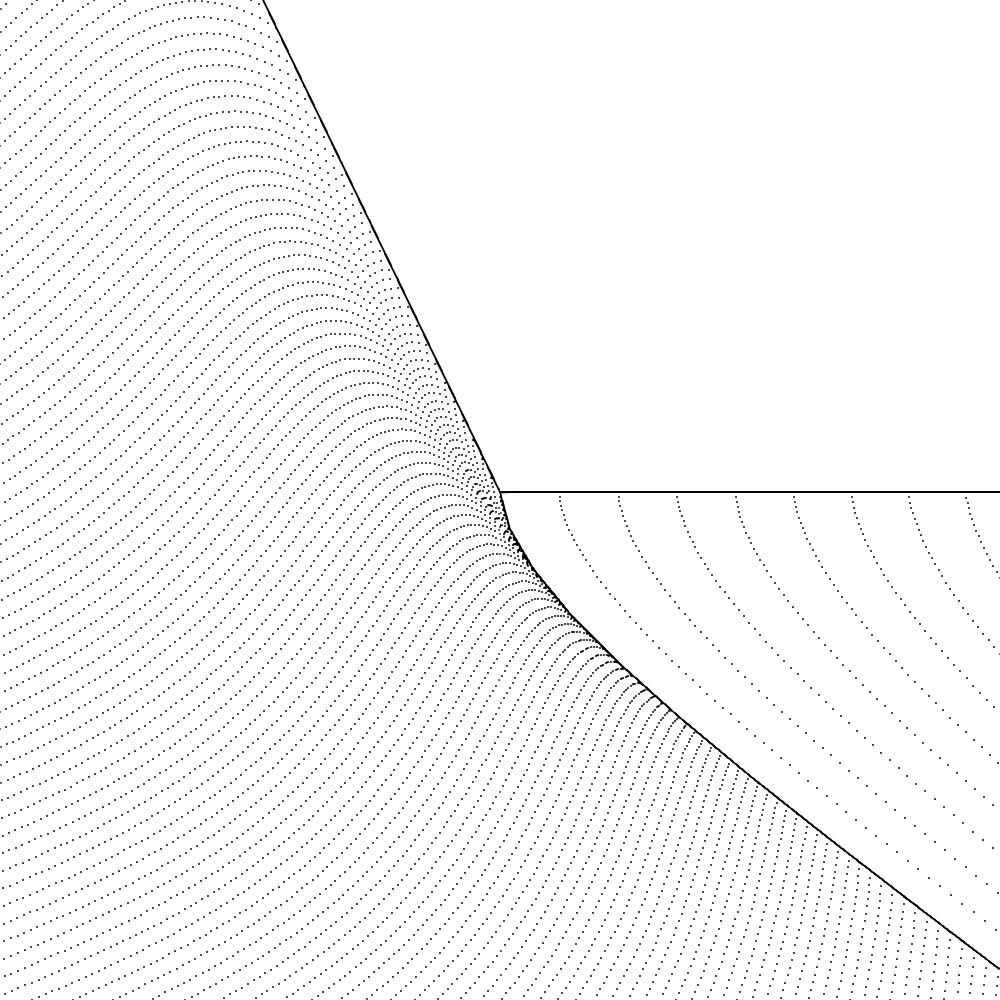}}}
        \newline
        \caption{Mesh optimization without (top row) and with (bottom row) the proposed mesh validity bounding strategy for a blade geometry. Original mesh shown on left, optimized mesh shown in middle, and zoom-in on red region shown on right. Zoom-in region shows node renderings (equispaced on the reference element) to visualize mesh inversion.}
        \label{fig:blade_mesh}
    \end{figure}

An example of this is shown in the $r$-adaptive mesh optimization of a two-dimensional turbine blade geometry using fourth-order quadrilateral elements. The optimization process with and without the proposed approach, which ensures mesh validity either at discrete locations (mesh nodes and quadrature points) or throughout the entire element, respectively, is shown in \cref{fig:blade_mesh}. Without the proposed approach, the optimization process yielded an inverted mesh element in the top right of the domain, which can be seen by twisting-like behavior in the visualization of nodes in the reference space of the element. With the proposed approach, this behavior was not observed, and the optimized mesh is guaranteed to be valid.

\subsection{Continuously bounds-preserving limiting}
The proposed approach was then implemented in the context of continuously bounds-preserving limiting for high-order DG schemes. We first look at a representative example of the interpolation of a step function, where the presence of Gibbs phenomena in polynomial interpolation results in overshoots and undershoots, a commonly encountered problem in the approximation of hyperbolic conservation laws. The metric of interest here is the ability of the approach in predicting a lower/upper bound on the solution polynomial, which can be compared to the true minima/maxima (obtained via brute force sampling approximations) and methods such as taking the minimum/maximum Bernstein coefficients. We consider an interpolation of a step function $u_0(x)$ onto the Gauss--Lobatto nodes, given as 
\begin{equation}
    u_0(x) = \begin{cases}
        -0.5 &\text{ if } x < 0,\\
        \hphantom{-.5}0 &\text{ if } x = 0,\\
        \hphantom{-}0.5 &\text{ if }x > 0.
    \end{cases}
\end{equation}
Here, the extrema of the solution polynomial $u_h(x)$ exceed $\pm 0.5$ due to overshoots/undershoots, with their magnitude depending on the order of the interpolating polynomial. A comparison of the predicted extrema as obtained by an exact (brute force sampling) approximation, the proposed approach (with $M = N$), and a Bernstein-based approach is shown in \cref{tab:step_bounds_error} for a variety of approximation orders. It can be seen that the proposed approach tightly bounds the true extrema, with relative errors generally less than $10\%$ that only mildly increased with increasing approximation order. In contrast, the Bernstein approximation yielded much larger errors, often $10-100$ times larger than the proposed approach, which grew quickly with increasing approximation order. These results indicate that the proposed method may yield much tighter bounds than approaches based on convex hull properties such as Bernstein representations. 

    \begin{figure}[htbp!] 
        \centering
        \begin{tabular*}{0.9\linewidth}{@{\extracolsep{\fill}} |r | ccccc | }
        \hline
         & $\mathbb P_3$ & $\mathbb P_4$ & $\mathbb P_5$ & $\mathbb P_6$  & $\mathbb P_7$ \\ 
        \hline
        Exact                  & $\pm\ 0.6286$ & $\pm\ 0.5342$ & $\pm\ 0.6368$ & $\pm\ 0.5340$ & $\pm\ 0.6389$\\
        Bernstein              & $\pm\ 1.1967$ & $\pm\ 0.7116$ & $\pm\ 3.0758$ & $\pm\ 1.1040$ & $\pm\ 8.8450$ \\
        Present work ($M = N$) & $\pm\ 0.6530$ & $\pm\ 0.5867$ & $\pm\ 0.6780$ & $\pm\ 0.6236$ & $\pm\ 0.7080$ \\
        \hline
        Error reduction & $-95.7\%$ & $-70.4\%$ & $-98.3\%$ & $-84.3\%$ & $-99.2\%$\\
        \hline
        \end{tabular*}
        \captionof{table}{\label{tab:step_bounds_error} Predicted extrema for a polynomial interpolating a step function $[-0.5, 0.5]$ on the Gauss--Lobatto nodes at varying approximation orders. Error reduction in the bounds between the present approach (with $M = N$) and Bernstein approach shown on bottom.}
    \end{figure}

This approach was then implemented in time-dependent simulations of the linear transport equation using the limiter described in \cref{ssec:limiter}. The example of the solid body rotation problem of \citet{LeVeque1996} was used, where the domain is set as the unit square $[0,1]^2$ with periodic boundary conditions, and the initial conditions are given as
\begin{equation}
    u_0(\mathbf{x}) = \begin{cases}
        1, &\text{if } (x - 0.5)^2 + (y - 0.75)^2 \leq 0.15^2 \\
        & \quad \text{ and } x,y \notin [0.475,0.525]\times[0.6, 0.85], \\
        0.25\left (1 + \cos \left (\frac{\pi}{0.15} \sqrt{(x - 0.25)^2 + (y - 0.5)^2}\right ) \right ), &\text{if } (x - 0.25)^2 + (y - 0.5)^2 \leq 0.15^2, \\
        1 - \sqrt{(x - 0.5)^2 + (y - 0.25)^2}/0.15, &\text{if } (x - 0.5)^2 + (y - 0.25)^2 \leq 0.15^2, \\
        0, &\text{else}.
    \end{cases} 
\end{equation}
The advection velocity field was set as  
\begin{equation}
    c(\mathbf{x}) = [-2 \pi (y - 0.5),\ 2 \pi (x - 0.5)]^T,
\end{equation}
which induced a counterclockwise rotation with constant angular velocity about the domain center, completing one full revolution per unit time. These initial conditions include solution profiles of varying continuity: a $C^\infty$ cosinusoidal hump, a $C^0$ sharp cone, and a discontinuous notched cylinder. These features pose challenges for high-order approximations, potentially leading to violations of the maximum principle. For this problem, a global maximum principle of the form $u_h(\mathbf{x},t) \in [0,1]$ was enforced continuously using the proposed approach. 

    \begin{figure}[htbp!]
        \centering
        \subfloat[$N_e = 32^2$]{
        \adjustbox{width=0.33\linewidth, valign=b}{\includegraphics[]{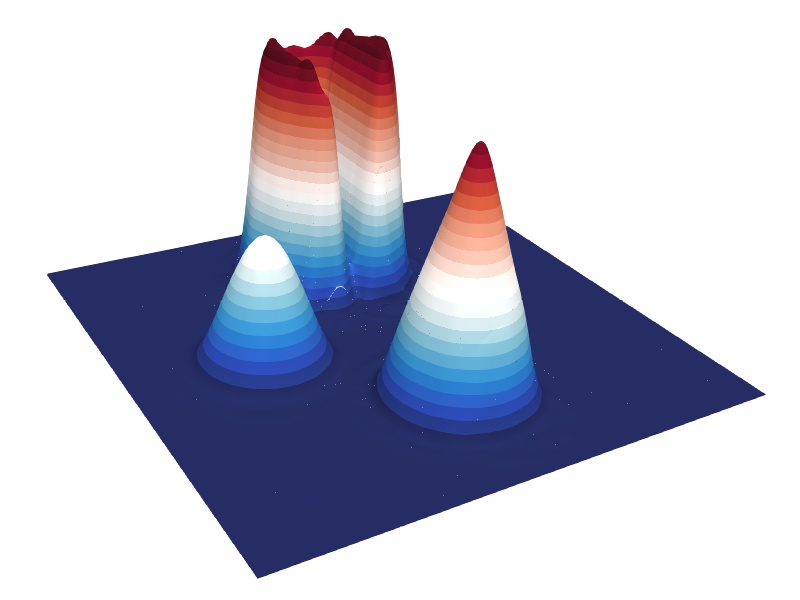}}}
        \subfloat[$N_e = 64^2$]{
        \adjustbox{width=0.33\linewidth, valign=b}{\includegraphics[]{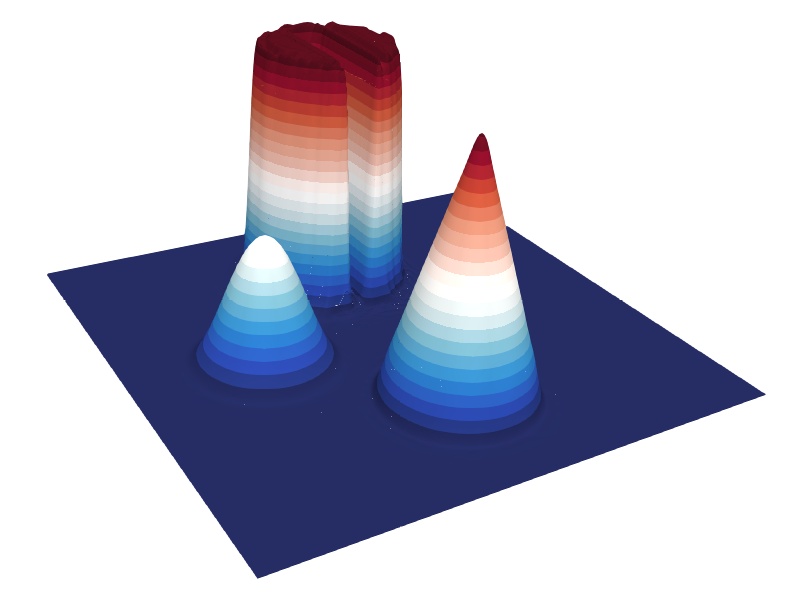}}}
        \subfloat[$N_e = 128^2$]{
        \adjustbox{width=0.33\linewidth, valign=b}{\includegraphics[]{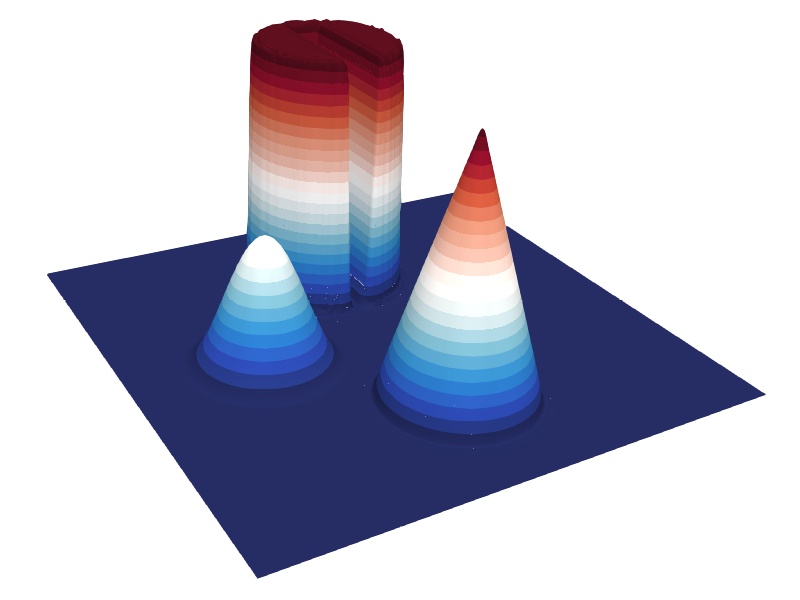}}}
        \caption{Solution contours for the solid body rotation problem at $t=1$ as computed by a $\mathbb P_3$ DG approximation with global maximum principle bounds $u_h(\mathbf{x},t) \in [0,1]$ on a varying number of quadrilateral mesh elements $N_e$. }
        \label{fig:sbr_profiles}
    \end{figure}

The solution as computed with a $\mathbb P_3$ DG approximation using uniform meshes with a varying number of quadrilateral elements (denoted by $N_e$) is shown in \cref{fig:sbr_profiles} at the final time $t=1$. It can be seen that with increasing resolution, the numerical diffusion around the notched cylinder decreased, such that the initial profile was well-recovered at the highest resolution. Furthermore, the solution did not exceed the bounds of the initial conditions. This was quantitatively verified and presented in \cref{tab:sbr_minmax}, which shows the minimum and maximum values of $u_h(\mathbf{x})$, computed via a brute force sampling approach, at the final time. The proposed approach ensured that the solution remained within the bounds $u_h(\mathbf{x},t) \in [0,1]$, with the highest resolution enforcing bounds to essentially machine precision levels. This convergence of the extrema to the exact bounds with respect to the increasing mesh resolution can be attributed to the lower numerical dissipation at higher resolution levels. 

    \begin{figure}[htbp!] 
        \centering
        \begin{tabular*}{0.4\linewidth}{@{\extracolsep{\fill}} |c | c | c | }
        \hline
         $N_e$ & $\underset{\mathbf{x}}{\min}\ u_h(\mathbf{x}, 1)$ & $\underset{\mathbf{x}}{\max}\ u_h(\mathbf{x}, 1)$ \\ 
        \hline
        $16^2$ & $9.08592 \times 10^{-7}$ & $0.97632323$\\
        $32^2$ & $1.56217 \times 10^{-9}$ & $0.98881219$\\
        $64^2$ & $8.34347 \times 10^{-13}$ & $0.99984185$\\
        $128^2$ & $7.86022 \times 10^{-14}$ & $0.99999856$\\
        \hline
        \end{tabular*}
        \captionof{table}{\label{tab:sbr_minmax} Minimum/maximum values of the solution (computed via brute force sampling) for the solid body rotation problem  at $t=1$ as computed by a $\mathbb P_3$ DG approximation with global maximum principle bounds $u_h(\mathbf{x},t) \in [0,1]$ on a varying number of quadrilateral mesh elements $N_e$.}
    \end{figure}

\section{Conclusions}\label{sec:conclusion}
In this work, we introduced a novel approach to bounding extrema in high-order polynomial approximations in finite element methods. The approach relies on precomputed piece-wise linear bounding boxes for polynomial basis functions stemming from a constrained optimization problem, enabling accurate and efficient local bounds for any polynomial formed from these bases. The accuracy of these bounds can be further improved through a simple basis transformation, and the method is applicable to arbitrary element types and approximation orders while remaining computationally efficient for on-the-fly evaluation. Some applications of the proposed approach are shown in mesh validity checks and optimization for high-order curved meshes, ensuring positivity of the element Jacobian determinant, as well as in continuously bounds-preserving limiters for hyperbolic systems, where it enforces maximum principle constraints across the entire solution polynomial. Furthermore, comparisons to traditional approaches relying on convex hull properties such as Bernstein polynomials show significantly tighter bounds. These results highlight the potential of the approach, which is applicable to not only finite element methods but also various problems in scientific computing ranging from computer graphics to collision detection.

\section*{Acknowledgements}
\label{sec:ack}
This work was  performed under the auspices of the U.S. Department of Energy by Lawrence Livermore National Laboratory under contract DE--AC52--07NA27344. Release number LLNL--JRNL--2004870--DRAFT.
\bibliographystyle{unsrtnat}
\bibliography{reference}

\clearpage
\begin{appendices}
\crefalias{section}{appendix}
\section{Tensor-product optimizations for higher dimensions}\label{app:higherdim}

For tensor-product bases in higher dimensions, computation of the bounds and corresponding $L^2$ projections simplifies to solving a sequence of one-dimensional problems along each coordinate direction. For example, the two-dimensional tensor-product solution, 
\begin{equation}
    u_h(x,y) = \sum_{j=1}^N\sum_{i=1}^N u_{ij} \phi_i(x) \phi_j(y) 
\end{equation}
can be equivalent represented in terms of a polynomial along one direction (e.g. $y$) for which its coefficients $w_j$ vary with respect to the other direction (e.g., $x$) as
\begin{subequations}\label{eq:2D_bounds}
    \begin{align}
    u(x,y) &= \sum_{j=1}^N\underbrace{\bigg(\sum_{i=1}^N u_{ij} \phi_i(x)\bigg)}_{w_j(x)} \phi_j(y),\\
    u(x=\eta_k,y) &= \sum_{j=1}^N w_j(x=\eta_k) \phi_j(y).\label{eq:2D_bounds_2}
    \end{align}
\end{subequations}
In \eqref{eq:2D_bounds_2}, the minimum and maximum bounds on $w_j(x=\eta_k)$ can be computed using \eqref{eq:1D_bounds_offset} at the $k=1\dots M$ control points for each $j \in [1,N]$. Then, bounding \eqref{eq:2D_bounds_2} entails repeating the one-dimensional bounding procedure:
\begin{subequations}\label{eq:2D_bounds_final}
    \begin{align}
    q^-(\eta_k, \eta_l) = u_{k,LO}(\eta_l) + \sum_{j=1}^N \min(w^{'-}_j(\eta_k) q_j^-(\eta_l),w^{'+}_j(\eta_k) q_j^-(\eta_l)),\\
    q^+(\eta_k, \eta_l) = u_{k,LO}(\eta_l) + \sum_{j=1}^N \max(w^{'-}_j(\eta_k) q_j^-(\eta_l),w^{'+}_j(\eta_k) q_j^-(\eta_l)).
    \end{align}
\end{subequations}
Here, $u_{k,LO}$ represents the $\mathbb P_1$ fit of $u(x=\eta_k,y)$ and $w^{'-}_j$ and $w^{'+}_j$ denote the minimum/maximum bounds, respectively, on the coefficients $w_j'$ of the high-order fluctuations (computed after offsetting the linear fit). 
The tensor-product simplification results in the total computational cost of $\mathcal O (N^{d}M + NM^{d})$ for a $d$-dimensional function (as opposed to $\mathcal O(N^dM^d)$). An implementation of these optimizations is presented in \citep{gslib-github}. Note that this simplification does not extend to non tensor-product elements (e.g., simplices) where the computational cost of bounding the polynomial then scales directly as the product of the total number of basis functions and control nodes.

\section{Bounding box examples}\label{app:bounds}

We present here some tabulated examples of the optimized bounding boxes for Gauss--Lobatto interpolating basis functions at varying approximation orders and number of control nodes from $M = N$ to $M = N+2$.

\subsection{$M = N$}
\begin{figure}[htbp!]  
\centering
\begin{tabular*}{0.9\linewidth}{@{\extracolsep{\fill}} |r | rrr | }
\hline
& $\phi_1$ & $\phi_2$ & $\phi_3$\\
\hline
$L^1$ & $0.8752491$ & $-0.0010021$ & $-0.1247509$\\
$L^2$ & $-0.1252514$ & $0.9999990$ & $-0.1252514$\\
$L^3$ & $-0.1247509$ & $-0.0010021$ & $0.8752491$\\
\hline
$U^1$ & $1.0005016$ & $0.2495008$ & $0.0005016$\\
$U^2$ & $0.0000010$ & $1.2505018$ & $0.0000010$\\
$U^3$ & $0.0005016$ & $0.2495008$ & $1.0005016$\\
\hline
\end{tabular*}
\captionof{table}{Lower and upper bounding box control node values for a $\mathbb P_2$ Gauss--Lobatto interpolating polynomial with $M = 3$. Control nodes $\boldsymbol{\eta} = [-1.0, 0.0, 1.0]$.}
\end{figure}

\begin{figure}[htbp!]  
\centering
\begin{tabular*}{0.9\linewidth}{@{\extracolsep{\fill}} |r | rrrr | }
\hline
& $\phi_1$ & $\phi_2$ & $\phi_3$ & $\phi_4$\\
\hline
$L^1$ & $0.8948521$ & $-0.0014658$ & $-0.3252945$ & $-0.0000031$\\
$L^2$ & $-0.1811453$ & $0.9993988$ & $-0.0794586$ & $-0.0688717$\\
$L^3$ & $-0.0688717$ & $-0.0794586$ & $0.9993988$ & $-0.1811453$\\
$L^4$ & $-0.0000031$ & $-0.3252945$ & $-0.0014658$ & $0.8948521$\\
\hline
$U^1$ & $1.0008087$ & $0.2667580$ & $0.0003650$ & $0.1253381$\\
$U^2$ & $0.0128333$ & $1.2677409$ & $-0.0169454$ & $0.0118659$\\
$U^3$ & $0.0118659$ & $-0.0169454$ & $1.2677409$ & $0.0128333$\\
$U^4$ & $0.1253381$ & $0.0003650$ & $0.2667580$ & $1.0008087$\\
\hline
\end{tabular*}
\captionof{table}{Lower and upper bounding box control node values for a $\mathbb P_3$ Gauss--Lobatto interpolating polynomial with $M = 4$. Control nodes $\boldsymbol{\eta} = [-1.0, -0.4626417, 0.4626417, 1.0]$.}
\end{figure}

\begin{figure}[htbp!]  
\centering
\begin{tabular*}{0.9\linewidth}{@{\extracolsep{\fill}} |r | rrrrr | }
\hline
& $\phi_1$ & $\phi_2$ & $\phi_3$ & $\phi_4$ & $\phi_5$\\
\hline
$L^1$ & $0.9026724$ & $-0.0023273$ & $-0.3515760$ & $-0.0023273$ & $-0.0738878$\\
$L^2$ & $-0.2002313$ & $0.9984396$ & $-0.0741776$ & $-0.2067862$ & $-0.0057559$\\
$L^3$ & $-0.0478948$ & $-0.1949668$ & $0.9999989$ & $-0.1949668$ & $-0.0478948$\\
$L^4$ & $-0.0057559$ & $-0.2067862$ & $-0.0741776$ & $0.9984396$ & $-0.2002313$\\
$L^5$ & $-0.0738878$ & $-0.0023273$ & $-0.3515759$ & $-0.0023273$ & $0.9026724$\\
\hline
$U^1$ & $1.0013235$ & $0.3227364$ & $0.0000027$ & $0.2049544$ & $0.0013235$\\
$U^2$ & $0.0210607$ & $1.2141881$ & $0.0060579$ & $0.0199534$ & $0.0909730$\\
$U^3$ & $0.0438391$ & $0.0000010$ & $1.3036861$ & $0.0000010$ & $0.0438391$\\
$U^4$ & $0.0909730$ & $0.0199534$ & $0.0060579$ & $1.2141881$ & $0.0210607$\\
$U^5$ & $0.0013235$ & $0.2049544$ & $0.0000027$ & $0.3227364$ & $1.0013235$\\
\hline
\end{tabular*}
\captionof{table}{Lower and upper bounding box control node values for a $\mathbb P_4$ Gauss--Lobatto interpolating polynomial with $M = 5$. Control nodes $\boldsymbol{\eta} = [-1.0, -0.6708529, 0.0, 0.6708529, 1.0]$.}
\end{figure}

\begin{figure}[htbp!]  
\centering
\begin{tabular*}{0.9\linewidth}{@{\extracolsep{\fill}} |r | rrrrrr | }
\hline
& $\phi_1$ & $\phi_2$ & $\phi_3$ & $\phi_4$ & $\phi_5$ & $\phi_6$\\
\hline
$L^1$ & $0.8513378$ & $-0.0028296$ & $-0.3920037$ & $-0.0021114$ & $-0.1421587$ & $-0.0000028$\\
$L^2$ & $-0.2170296$ & $0.9995489$ & $-0.0323278$ & $-0.2548223$ & $-0.0007090$ & $-0.0607456$\\
$L^3$ & $-0.0224756$ & $-0.2453042$ & $0.9999529$ & $-0.1675032$ & $-0.1511019$ & $-0.0189264$\\
$L^4$ & $-0.0189264$ & $-0.1511019$ & $-0.1675032$ & $0.9999529$ & $-0.2453042$ & $-0.0224756$\\
$L^5$ & $-0.0607456$ & $-0.0007090$ & $-0.2548223$ & $-0.0323278$ & $0.9995489$ & $-0.2170296$\\
$L^6$ & $-0.0000028$ & $-0.1421587$ & $-0.0021114$ & $-0.3920037$ & $-0.0028296$ & $0.8513378$\\
\hline
$U^1$ & $1.0021237$ & $0.4395637$ & $0.0004805$ & $0.2277843$ & $0.0007634$ & $0.0535520$\\
$U^2$ & $-0.0106906$ & $1.1844705$ & $0.1323575$ & $0.0008303$ & $0.1665040$ & $0.0012256$\\
$U^3$ & $0.0633344$ & $-0.0049788$ & $1.2671846$ & $0.0051295$ & $0.0586046$ & $0.0687812$\\
$U^4$ & $0.0687812$ & $0.0586046$ & $0.0051295$ & $1.2671846$ & $-0.0049788$ & $0.0633344$\\
$U^5$ & $0.0012256$ & $0.1665040$ & $0.0008303$ & $0.1323575$ & $1.1844705$ & $-0.0106906$\\
$U^6$ & $0.0535520$ & $0.0007634$ & $0.2277843$ & $0.0004805$ & $0.4395637$ & $1.0021237$\\
\hline
\end{tabular*}
\captionof{table}{Lower and upper bounding box control node values for a $\mathbb P_5$ Gauss--Lobatto interpolating polynomial with $M = 6$. Control nodes $\boldsymbol{\eta} = [-1.0, -0.7589109, -0.2823207, 0.2823207, 0.7589109, 1.0]$.}
\end{figure}

\begin{figure}[htbp!]  
\centering
\begin{tabular*}{0.9\linewidth}{@{\extracolsep{\fill}} |r | rrrrrrr | }
\hline
& $\phi_1$ & $\phi_2$ & $\phi_3$ & $\phi_4$ & $\phi_5$ & $\phi_6$ & $\phi_7$\\
\hline
$L^1$ & $0.8042527$ & $-0.0037149$ & $-0.4107312$ & $-0.0024327$ & $-0.1579100$ & $-0.0037149$ & $-0.0368755$\\
$L^2$ & $-0.2275404$ & $0.9950029$ & $-0.0004258$ & $-0.2680649$ & $0.0105489$ & $-0.1254623$ & $0.0031652$\\
$L^3$ & $-0.0133993$ & $-0.2691138$ & $0.9994716$ & $-0.1608572$ & $-0.1963252$ & $-0.0309949$ & $-0.0486212$\\
$L^4$ & $-0.0312508$ & $-0.1278783$ & $-0.2224863$ & $0.9999990$ & $-0.2224863$ & $-0.1278783$ & $-0.0312508$\\
$L^5$ & $-0.0486212$ & $-0.0309949$ & $-0.1963252$ & $-0.1608572$ & $0.9994716$ & $-0.2691138$ & $-0.0133993$\\
$L^6$ & $0.0031652$ & $-0.1254623$ & $0.0105489$ & $-0.2680649$ & $-0.0004258$ & $0.9950029$ & $-0.2275404$\\
$L^7$ & $-0.0368755$ & $-0.0037149$ & $-0.1579100$ & $-0.0024327$ & $-0.4107312$ & $-0.0037149$ & $0.8042527$\\
\hline
$U^1$ & $1.0028231$ & $0.4976582$ & $0.0006166$ & $0.2356049$ & $0.0006166$ & $0.0993084$ & $0.0028231$\\
$U^2$ & $-0.0339076$ & $1.1872082$ & $0.1688978$ & $-0.0154549$ & $0.1829468$ & $-0.0073889$ & $0.0489994$\\
$U^3$ & $0.0744049$ & $0.0173077$ & $1.2513780$ & $-0.0170501$ & $0.0580635$ & $0.1376125$ & $0.0123804$\\
$U^4$ & $0.0576242$ & $0.0889146$ & $0.0000010$ & $1.3036627$ & $0.0000010$ & $0.0889146$ & $0.0576242$\\
$U^5$ & $0.0123804$ & $0.1376125$ & $0.0580635$ & $-0.0170501$ & $1.2513780$ & $0.0173077$ & $0.0744049$\\
$U^6$ & $0.0489994$ & $-0.0073889$ & $0.1829468$ & $-0.0154549$ & $0.1688978$ & $1.1872082$ & $-0.0339076$\\
$U^7$ & $0.0028231$ & $0.0993084$ & $0.0006166$ & $0.2356049$ & $0.0006166$ & $0.4976582$ & $1.0028231$\\
\hline
\end{tabular*}
\captionof{table}{Lower and upper bounding box control node values for a $\mathbb P_6$ Gauss--Lobatto interpolating polynomial with $M = 7$. Control nodes $\boldsymbol{\eta} = [-1.0, -0.815025, -0.476498, 0.0, 0.476498, 0.815025, 1.0]$.}
\end{figure}

\newpage
\subsection{$M = N + 1$}
\begin{figure}[htbp!]  
\centering
\begin{tabular*}{0.9\linewidth}{@{\extracolsep{\fill}} |r | rrr | }
\hline
& $\phi_1$ & $\phi_2$ & $\phi_3$\\
\hline
$L^1$ & $0.9451077$ & $-0.0006684$ & $-0.0548923$\\
$L^2$ & $0.1671654$ & $0.8882205$ & $-0.1671675$\\
$L^3$ & $-0.1671675$ & $0.8882205$ & $0.1671654$\\
$L^4$ & $-0.0548923$ & $-0.0006684$ & $0.9451077$\\
\hline
$U^1$ & $1.0003347$ & $0.1097835$ & $0.0003347$\\
$U^2$ & $0.2230567$ & $1.0000011$ & $-0.1112761$\\
$U^3$ & $-0.1112761$ & $1.0000010$ & $0.2230567$\\
$U^4$ & $0.0003347$ & $0.1097835$ & $1.0003347$\\
\hline
\end{tabular*}
\captionof{table}{Lower and upper bounding box control node values for a $\mathbb P_2$ Gauss--Lobatto interpolating polynomial with $M = 4$. Control nodes $\boldsymbol{\eta} = [-1.0, -0.3343328, 0.3343328, 1.0]$.}
\end{figure}

\begin{figure}[htbp!]  
\centering
\begin{tabular*}{0.9\linewidth}{@{\extracolsep{\fill}} |r | rrrr | }
\hline
& $\phi_1$ & $\phi_2$ & $\phi_3$ & $\phi_4$\\
\hline
$L^1$ & $0.9133175$ & $-0.0013835$ & $-0.1691465$ & $-0.0001780$\\
$L^2$ & $-0.0060187$ & $0.9657676$ & $-0.2000494$ & $0.0314012$\\
$L^3$ & $-0.1895477$ & $0.6236876$ & $0.6236876$ & $-0.1895477$\\
$L^4$ & $0.0314012$ & $-0.2000494$ & $0.9657676$ & $-0.0060187$\\
$L^5$ & $-0.0001780$ & $-0.1691465$ & $-0.0013835$ & $0.9133175$\\
\hline
$U^1$ & $1.0007523$ & $0.1722848$ & $0.0007549$ & $0.0632003$\\
$U^2$ & $0.1115588$ & $1.1738419$ & $-0.1087275$ & $0.0531246$\\
$U^3$ & $-0.1244389$ & $0.7053872$ & $0.7053872$ & $-0.1244389$\\
$U^4$ & $0.0531246$ & $-0.1087275$ & $1.1738419$ & $0.1115588$\\
$U^5$ & $0.0632003$ & $0.0007549$ & $0.1722848$ & $1.0007523$\\
\hline
\end{tabular*}
\captionof{table}{Lower and upper bounding box control node values for a $\mathbb P_3$ Gauss--Lobatto interpolating polynomial with $M = 5$. Control nodes $\boldsymbol{\eta} = [-1.0, -0.5606852, 0.0, 0.5606852, 1.0]$.}
\end{figure}

\begin{figure}[htbp!]  
\centering
\begin{tabular*}{0.9\linewidth}{@{\extracolsep{\fill}} |r | rrrrr | }
\hline
& $\phi_1$ & $\phi_2$ & $\phi_3$ & $\phi_4$ & $\phi_5$\\
\hline
$L^1$ & $0.9205495$ & $-0.0022944$ & $-0.2394212$ & $-0.0022944$ & $-0.0540298$\\
$L^2$ & $-0.1007104$ & $0.9817825$ & $-0.1552628$ & $-0.0006023$ & $-0.0190016$\\
$L^3$ & $-0.1602338$ & $0.2756817$ & $0.9011641$ & $-0.3140172$ & $0.0500969$\\
$L^4$ & $0.0500969$ & $-0.3140172$ & $0.9011641$ & $0.2756817$ & $-0.1602338$\\
$L^5$ & $-0.0190016$ & $-0.0006023$ & $-0.1552628$ & $0.9817825$ & $-0.1007104$\\
$L^6$ & $-0.0540298$ & $-0.0022944$ & $-0.2394212$ & $-0.0022944$ & $0.9205495$\\
\hline
$U^1$ & $1.0012947$ & $0.1944317$ & $0.0007171$ & $0.1535178$ & $0.0012947$\\
$U^2$ & $0.0783273$ & $1.2354373$ & $-0.0857886$ & $0.0557246$ & $0.0087225$\\
$U^3$ & $-0.0699741$ & $0.3570428$ & $1.0007138$ & $-0.1656079$ & $0.0947372$\\
$U^4$ & $0.0947372$ & $-0.1656079$ & $1.0007138$ & $0.3570428$ & $-0.0699741$\\
$U^5$ & $0.0087225$ & $0.0557246$ & $-0.0857886$ & $1.2354373$ & $0.0783273$\\
$U^6$ & $0.0012947$ & $0.1535178$ & $0.0007171$ & $0.1944317$ & $1.0012947$\\
\hline
\end{tabular*}
\captionof{table}{Lower and upper bounding box control node values for a $\mathbb P_4$ Gauss--Lobatto interpolating polynomial with $M = 6$. Control nodes $\boldsymbol{\eta} = [-1.0, -0.7088516, -0.1740483, 0.1740483, 0.7088516, 1.0]$.}
\end{figure}

\begin{figure}[htbp!]  
\centering
\begin{tabular*}{0.9\linewidth}{@{\extracolsep{\fill}} |r | rrrrrr | }
\hline
& $\phi_1$ & $\phi_2$ & $\phi_3$ & $\phi_4$ & $\phi_5$ & $\phi_6$\\
\hline
$L^1$ & $0.8814769$ & $-0.0029636$ & $-0.3117680$ & $-0.0017280$ & $-0.1154057$ & $-0.0002392$\\
$L^2$ & $-0.1585389$ & $0.9968812$ & $-0.0992253$ & $-0.0897185$ & $-0.0204307$ & $-0.0230939$\\
$L^3$ & $-0.1045286$ & $0.0675661$ & $0.9617801$ & $-0.2884231$ & $0.0557476$ & $-0.0484811$\\
$L^4$ & $0.0624990$ & $-0.3492784$ & $0.6306984$ & $0.6306984$ & $-0.3492784$ & $0.0624990$\\
$L^5$ & $-0.0484811$ & $0.0557476$ & $-0.2884231$ & $0.9617801$ & $0.0675661$ & $-0.1045286$\\
$L^6$ & $-0.0230939$ & $-0.0204307$ & $-0.0897185$ & $-0.0992253$ & $0.9968812$ & $-0.1585389$\\
$L^7$ & $-0.0002392$ & $-0.1154057$ & $-0.0017280$ & $-0.3117680$ & $-0.0029636$ & $0.8814769$\\
\hline
$U^1$ & $1.0021031$ & $0.2958899$ & $0.0011181$ & $0.1940650$ & $0.0006478$ & $0.0433088$\\
$U^2$ & $0.0302645$ & $1.2313505$ & $-0.0380837$ & $0.0317172$ & $0.0624364$ & $0.0084858$\\
$U^3$ & $-0.0307161$ & $0.1666680$ & $1.1549303$ & $-0.1222708$ & $0.1379809$ & $-0.0183029$\\
$U^4$ & $0.1220502$ & $-0.1946476$ & $0.7419176$ & $0.7419176$ & $-0.1946476$ & $0.1220502$\\
$U^5$ & $-0.0183029$ & $0.1379809$ & $-0.1222708$ & $1.1549303$ & $0.1666680$ & $-0.0307161$\\
$U^6$ & $0.0084858$ & $0.0624364$ & $0.0317172$ & $-0.0380837$ & $1.2313505$ & $0.0302645$\\
$U^7$ & $0.0433088$ & $0.0006478$ & $0.1940650$ & $0.0011181$ & $0.2958899$ & $1.0021031$\\
\hline
\end{tabular*}
\captionof{table}{Lower and upper bounding box control node values for a $\mathbb P_5$ Gauss--Lobatto interpolating polynomial with $M = 7$. Control nodes $\boldsymbol{\eta} = [-1.0, -0.7809033, -0.3683177, 0.0, 0.3683177, 0.7809033, 1.0]$.}
\end{figure}

\begin{figure}[htbp!]  
\centering
\begin{tabular*}{0.9\linewidth}{@{\extracolsep{\fill}} |r | rrrrrrr | }
\hline
& $\phi_1$ & $\phi_2$ & $\phi_3$ & $\phi_4$ & $\phi_5$ & $\phi_6$ & $\phi_7$\\
\hline
$L^1$ & $0.8672816$ & $-0.0044051$ & $-0.3464291$ & $-0.0012834$ & $-0.1320094$ & $-0.0044051$ & $-0.0316918$\\
$L^2$ & $-0.1882615$ & $0.9994606$ & $-0.0715195$ & $-0.1477100$ & $-0.0133397$ & $-0.0720893$ & $-0.0024119$\\
$L^3$ & $-0.0701976$ & $-0.0705482$ & $0.9808790$ & $-0.2530301$ & $0.0389030$ & $-0.0729852$ & $0.0012625$\\
$L^4$ & $0.0435990$ & $-0.3213648$ & $0.3225332$ & $0.8711237$ & $-0.3511253$ & $0.0900655$ & $-0.0614205$\\
$L^5$ & $-0.0614205$ & $0.0900655$ & $-0.3511253$ & $0.8711237$ & $0.3225332$ & $-0.3213648$ & $0.0435990$\\
$L^6$ & $0.0012625$ & $-0.0729852$ & $0.0389030$ & $-0.2530301$ & $0.9808790$ & $-0.0705482$ & $-0.0701976$\\
$L^7$ & $-0.0024119$ & $-0.0720893$ & $-0.0133397$ & $-0.1477100$ & $-0.0715195$ & $0.9994606$ & $-0.1882615$\\
$L^8$ & $-0.0316918$ & $-0.0044051$ & $-0.1320094$ & $-0.0012834$ & $-0.3464291$ & $-0.0044051$ & $0.8672816$\\
\hline
$U^1$ & $1.0025747$ & $0.3876858$ & $0.0014470$ & $0.2088437$ & $0.0014470$ & $0.0842040$ & $0.0025747$\\
$U^2$ & $0.0122092$ & $1.1873133$ & $-0.0135598$ & $0.0199451$ & $0.1099858$ & $0.0068122$ & $0.0278756$\\
$U^3$ & $0.0074195$ & $0.1135189$ & $1.2483446$ & $-0.0924604$ & $0.1171454$ & $-0.0047314$ & $0.0280131$\\
$U^4$ & $0.1198678$ & $-0.1305006$ & $0.4334031$ & $1.0010918$ & $-0.1864547$ & $0.1717669$ & $-0.0307710$\\
$U^5$ & $-0.0307710$ & $0.1717669$ & $-0.1864547$ & $1.0010918$ & $0.4334031$ & $-0.1305006$ & $0.1198678$\\
$U^6$ & $0.0280131$ & $-0.0047314$ & $0.1171454$ & $-0.0924604$ & $1.2483446$ & $0.1135189$ & $0.0074195$\\
$U^7$ & $0.0278756$ & $0.0068122$ & $0.1099858$ & $0.0199451$ & $-0.0135598$ & $1.1873133$ & $0.0122092$\\
$U^8$ & $0.0025747$ & $0.0842040$ & $0.0014470$ & $0.2088437$ & $0.0014470$ & $0.3876858$ & $1.0025747$\\
\hline
\end{tabular*}
\captionof{table}{Lower and upper bounding box control node values for a $\mathbb P_6$ Gauss--Lobatto interpolating polynomial with $M = 8$. Control nodes $\boldsymbol{\eta} = [-1.0, -0.8350809, -0.5144829, -0.1380507, 0.1380507, 0.5144829, 0.8350809, 1.0]$.}
\end{figure}

\newpage
\subsection{$M = N + 2$}
\begin{figure}[htbp!]  
\centering
\begin{tabular*}{0.9\linewidth}{@{\extracolsep{\fill}} |r | rrr | }
\hline
& $\phi_1$ & $\phi_2$ & $\phi_3$\\
\hline
$L^1$ & $0.9689364$ & $-0.0005010$ & $-0.0310636$\\
$L^2$ & $0.3441870$ & $0.7494982$ & $-0.1563135$\\
$L^3$ & $-0.0313136$ & $0.9999989$ & $-0.0313136$\\
$L^4$ & $-0.1563135$ & $0.7494982$ & $0.3441870$\\
$L^5$ & $-0.0310636$ & $-0.0005010$ & $0.9689364$\\
\hline
$U^1$ & $1.0002510$ & $0.0621262$ & $0.0002510$\\
$U^2$ & $0.3755017$ & $0.8121255$ & $-0.1249988$\\
$U^3$ & $0.0000010$ & $1.0626262$ & $0.0000010$\\
$U^4$ & $-0.1249988$ & $0.8121255$ & $0.3755017$\\
$U^5$ & $0.0002510$ & $0.0621262$ & $1.0002510$\\
\hline
\end{tabular*}
\captionof{table}{Lower and upper bounding box control node values for a $\mathbb P_2$ Gauss--Lobatto interpolating polynomial with $M = 5$. Control nodes $\boldsymbol{\eta} = [-1.0, -0.5005005, 0.0, 0.5005005, 1.0]$.}
\end{figure}

\begin{figure}[htbp!]  
\centering
\begin{tabular*}{0.9\linewidth}{@{\extracolsep{\fill}} |r | rrrr | }
\hline
& $\phi_1$ & $\phi_2$ & $\phi_3$ & $\phi_4$\\
\hline
$L^1$ & $0.9293782$ & $-0.0012073$ & $-0.0763332$ & $-0.0002401$\\
$L^2$ & $0.1943347$ & $0.8698771$ & $-0.2587864$ & $0.0504246$\\
$L^3$ & $-0.1629644$ & $0.9291326$ & $0.2080413$ & $-0.0876780$\\
$L^4$ & $-0.0876780$ & $0.2080413$ & $0.9291326$ & $-0.1629644$\\
$L^5$ & $0.0504246$ & $-0.2587864$ & $0.8698771$ & $0.1943347$\\
$L^6$ & $-0.0002401$ & $-0.0763332$ & $-0.0012073$ & $0.9293782$\\
\hline
$U^1$ & $1.0006480$ & $0.1231346$ & $0.0006972$ & $0.0320353$\\
$U^2$ & $0.2486074$ & $0.9819215$ & $-0.1689090$ & $0.0728559$\\
$U^3$ & $-0.1004879$ & $1.0175813$ & $0.2358163$ & $-0.0571863$\\
$U^4$ & $-0.0571863$ & $0.2358163$ & $1.0175813$ & $-0.1004879$\\
$U^5$ & $0.0728559$ & $-0.1689090$ & $0.9819215$ & $0.2486074$\\
$U^6$ & $0.0320353$ & $0.0006972$ & $0.1231346$ & $1.0006480$\\
\hline
\end{tabular*}
\captionof{table}{Lower and upper bounding box control node values for a $\mathbb P_3$ Gauss--Lobatto interpolating polynomial with $M = 6$. Control nodes $\boldsymbol{\eta} = [-1.0, -0.6627409, -0.2704891, 0.2704891, 0.6627409, 1.0]$.}
\end{figure}

\begin{figure}[htbp!]  
\centering
\begin{tabular*}{0.9\linewidth}{@{\extracolsep{\fill}} |r | rrrrr | }
\hline
& $\phi_1$ & $\phi_2$ & $\phi_3$ & $\phi_4$ & $\phi_5$\\
\hline
$L^1$ & $0.9256847$ & $-0.0021771$ & $-0.1512635$ & $-0.0021771$ & $-0.0340882$\\
$L^2$ & $0.0332879$ & $0.9424694$ & $-0.2213716$ & $0.0614982$ & $-0.0321799$\\
$L^3$ & $-0.1970720$ & $0.7037976$ & $0.5453535$ & $-0.2638854$ & $0.0575199$\\
$L^4$ & $-0.0000277$ & $-0.0926919$ & $0.9999989$ & $-0.0926919$ & $-0.0000277$\\
$L^5$ & $0.0575199$ & $-0.2638854$ & $0.5453535$ & $0.7037976$ & $-0.1970720$\\
$L^6$ & $-0.0321799$ & $0.0614982$ & $-0.2213716$ & $0.9424694$ & $0.0332879$\\
$L^7$ & $-0.0340882$ & $-0.0021771$ & $-0.1512635$ & $-0.0021771$ & $0.9256847$\\
\hline
$U^1$ & $1.0012113$ & $0.1487162$ & $0.0011010$ & $0.0953204$ & $0.0012113$\\
$U^2$ & $0.1515358$ & $1.1429299$ & $-0.1355102$ & $0.0943803$ & $-0.0206778$\\
$U^3$ & $-0.1302360$ & $0.7846015$ & $0.6027643$ & $-0.1783198$ & $0.0886312$\\
$U^4$ & $0.0213470$ & $0.0000010$ & $1.1592583$ & $0.0000010$ & $0.0213470$\\
$U^5$ & $0.0886312$ & $-0.1783198$ & $0.6027643$ & $0.7846015$ & $-0.1302360$\\
$U^6$ & $-0.0206778$ & $0.0943803$ & $-0.1355102$ & $1.1429299$ & $0.1515358$\\
$U^7$ & $0.0012113$ & $0.0953204$ & $0.0011010$ & $0.1487162$ & $1.0012113$\\
\hline
\end{tabular*}
\captionof{table}{Lower and upper bounding box control node values for a $\mathbb P_4$ Gauss--Lobatto interpolating polynomial with $M = 7$. Control nodes $\boldsymbol{\eta} = [-1.0, -0.748869, -0.3902752, 0.0, 0.3902752, 0.748869, 1.0]$.}
\end{figure}

\begin{figure}[htbp!]  
\centering
\begin{tabular*}{0.9\linewidth}{@{\extracolsep{\fill}} |r | rrrrrr | }
\hline
& $\phi_1$ & $\phi_2$ & $\phi_3$ & $\phi_4$ & $\phi_5$ & $\phi_6$\\
\hline
$L^1$ & $0.9217631$ & $-0.0029380$ & $-0.2086252$ & $-0.0012131$ & $-0.0852034$ & $-0.0003771$\\
$L^2$ & $-0.0590266$ & $0.9696819$ & $-0.1812597$ & $0.0503179$ & $-0.0446515$ & $0.0039713$\\
$L^3$ & $-0.1835269$ & $0.4485174$ & $0.7816147$ & $-0.3214941$ & $0.1006510$ & $-0.0562234$\\
$L^4$ & $0.0354455$ & $-0.2631267$ & $0.9620203$ & $0.1327682$ & $-0.1330200$ & $0.0232467$\\
$L^5$ & $0.0232467$ & $-0.1330200$ & $0.1327682$ & $0.9620203$ & $-0.2631267$ & $0.0354455$\\
$L^6$ & $-0.0562234$ & $0.1006510$ & $-0.3214941$ & $0.7816147$ & $0.4485174$ & $-0.1835269$\\
$L^7$ & $0.0039713$ & $-0.0446515$ & $0.0503179$ & $-0.1812597$ & $0.9696819$ & $-0.0590266$\\
$L^8$ & $-0.0003771$ & $-0.0852034$ & $-0.0012131$ & $-0.2086252$ & $-0.0029380$ & $0.9217631$\\
\hline
$U^1$ & $1.0016678$ & $0.1777987$ & $0.0013446$ & $0.1391356$ & $0.0011777$ & $0.0320639$\\
$U^2$ & $0.1044963$ & $1.2121039$ & $-0.1067707$ & $0.0742302$ & $-0.0188173$ & $0.0168841$\\
$U^3$ & $-0.1032605$ & $0.5296652$ & $0.8683576$ & $-0.2021384$ & $0.1654181$ & $-0.0349335$\\
$U^4$ & $0.0758169$ & $-0.1195865$ & $1.1014727$ & $0.1696700$ & $-0.0682086$ & $0.0525290$\\
$U^5$ & $0.0525290$ & $-0.0682086$ & $0.1696700$ & $1.1014727$ & $-0.1195865$ & $0.0758169$\\
$U^6$ & $-0.0349335$ & $0.1654181$ & $-0.2021384$ & $0.8683576$ & $0.5296652$ & $-0.1032605$\\
$U^7$ & $0.0168841$ & $-0.0188173$ & $0.0742302$ & $-0.1067707$ & $1.2121039$ & $0.1044963$\\
$U^8$ & $0.0320639$ & $0.0011777$ & $0.1391356$ & $0.0013446$ & $0.1777987$ & $1.0016678$\\
\hline
\end{tabular*}
\captionof{table}{Lower and upper bounding box control node values for a $\mathbb P_5$ Gauss--Lobatto interpolating polynomial with $M = 8$. Control nodes $\boldsymbol{\eta} = [-1.0, -0.8130556, -0.4844534, -0.2002159, 0.2002159, 0.4844534, 0.8130556, 1.0]$.}
\end{figure}

\begin{figure}[htbp!]  
\centering
\begin{tabular*}{0.9\linewidth}{@{\extracolsep{\fill}} |r | rrrrrrr | }
\hline
& $\phi_1$ & $\phi_2$ & $\phi_3$ & $\phi_4$ & $\phi_5$ & $\phi_6$ & $\phi_7$\\
\hline
$L^1$ & $0.8820262$ & $-0.0043969$ & $-0.2829537$ & $-0.0018566$ & $-0.1144250$ & $-0.0043969$ & $-0.0275098$\\
$L^2$ & $-0.1316740$ & $0.9933112$ & $-0.1219835$ & $-0.0390320$ & $-0.0269245$ & $-0.0290033$ & $-0.0064344$\\
$L^3$ & $-0.1334867$ & $0.2069536$ & $0.9108360$ & $-0.3155131$ & $0.0852977$ & $-0.1012546$ & $0.0186195$\\
$L^4$ & $0.0599605$ & $-0.3253191$ & $0.8260690$ & $0.4037208$ & $-0.2764464$ & $0.0812513$ & $-0.0575199$\\
$L^5$ & $-0.0153225$ & $-0.0043037$ & $-0.1241529$ & $0.9999990$ & $-0.1241529$ & $-0.0043037$ & $-0.0153225$\\
$L^6$ & $-0.0575199$ & $0.0812513$ & $-0.2764464$ & $0.4037208$ & $0.8260690$ & $-0.3253191$ & $0.0599605$\\
$L^7$ & $0.0186195$ & $-0.1012546$ & $0.0852977$ & $-0.3155131$ & $0.9108360$ & $0.2069536$ & $-0.1334867$\\
$L^8$ & $-0.0064344$ & $-0.0290033$ & $-0.0269245$ & $-0.0390320$ & $-0.1219835$ & $0.9933112$ & $-0.1316740$\\
$L^9$ & $-0.0275098$ & $-0.0043969$ & $-0.1144250$ & $-0.0018566$ & $-0.2829537$ & $-0.0043969$ & $0.8820262$\\
\hline
$U^1$ & $1.0025552$ & $0.2731979$ & $0.0010892$ & $0.1817791$ & $0.0010892$ & $0.0729928$ & $0.0025552$\\
$U^2$ & $0.0453892$ & $1.2302202$ & $-0.0552581$ & $0.0437034$ & $0.0383789$ & $0.0171716$ & $0.0119274$\\
$U^3$ & $-0.0641939$ & $0.2868293$ & $1.0638067$ & $-0.1670870$ & $0.1621118$ & $-0.0502689$ & $0.0385482$\\
$U^4$ & $0.1066365$ & $-0.1927206$ & $0.9277952$ & $0.4688045$ & $-0.1547364$ & $0.1603258$ & $-0.0284066$\\
$U^5$ & $0.0025432$ & $0.0455342$ & $0.0002050$ & $1.1968939$ & $0.0002050$ & $0.0455342$ & $0.0025432$\\
$U^6$ & $-0.0284066$ & $0.1603258$ & $-0.1547364$ & $0.4688045$ & $0.9277952$ & $-0.1927206$ & $0.1066365$\\
$U^7$ & $0.0385482$ & $-0.0502689$ & $0.1621118$ & $-0.1670870$ & $1.0638067$ & $0.2868293$ & $-0.0641939$\\
$U^8$ & $0.0119274$ & $0.0171716$ & $0.0383789$ & $0.0437034$ & $-0.0552581$ & $1.2302202$ & $0.0453892$\\
$U^9$ & $0.0025552$ & $0.0729928$ & $0.0010892$ & $0.1817791$ & $0.0010892$ & $0.2731979$ & $1.0025552$\\
\hline
\end{tabular*}
\captionof{table}{Lower and upper bounding box control node values for a $\mathbb P_6$ Gauss--Lobatto interpolating polynomial with $M = 9$. Control nodes $\boldsymbol{\eta} = [-1.0, -0.8470722, -0.5666633, -0.3204544, 0.0, 0.3204544, 0.5666633, 0.8470722, 1.0]$.}
\end{figure}

\end{appendices}


\end{document}